\documentclass[12pt,twoside,a4paper]{amsart}
\usepackage[centertags]{amsmath}
\usepackage{amssymb,amscd,amsfonts,latexsym,a4wide,amsthm}

\advance\headheight by 3pt

\renewcommand{\subjclass}[1]{\thanks{\emph{2000 Mathematics Subject Classification:}~#1}}
\renewcommand{\keywords}[1]{\thanks{\emph{Keywords and Phrases:}~#1}}

\newtheorem{thm}{Theorem}[section]
\newtheorem{cor}[thm]{Corollary}
\newtheorem{lemma}[thm]{Lemma}
\newtheorem{prop}[thm]{Proposition}

\theoremstyle{definition}

\newtheorem{prg}{}[section]

\theoremstyle{remark}

\numberwithin{equation}{section}

\newcommand{\thmref}[1]{Theorem~\ref{#1}}
\newcommand{\secref}[1]{Section~\ref{#1}}
\newcommand{\lemref}[1]{Lemma~\ref{#1}}
\newcommand{\propref}[1]{Proposition~\ref{#1}}

\newcommand{\formref}[1]{(\ref{#1})}
\newcommand{\prgref}[1]{{\bf \S\ref{#1}}}

\def\rep#1{\bysame}
\def\bib[#1]#2<#3>#4|#5(#6)#7-#8.{\bibitem{#1}
{\sc #2},\ #3,\ {\it #4},\ {\bf #5}\ (#6)\ #7--#8.}
\def\tbib[#1]#2<#3>#4.{\bibitem{#1} {\sc #2},\ #3,\ #4.}
\def\bibt[#1]#2<#3>#4|#5(#6)#7-#8{\bibitem{#1}
{\sc #2},\ #3,\ in\ {\it #4},\ ed.\ {#5},\ (#6)\ #7--#8.}
\def\bibit[#1]#2<#3>{\bibitem{#1} {\sc#2},\ #3}
\def\prebib[#1]#2<#3>{\bibitem{#1} {\sc #2},\ {\it #3}, Preprint.}
\def\toabib[#1]#2<#3>(#4){\bibitem{#1} {\sc #2},\ #3,\ Accepted for
publication in {\it #4}}

\newcommand{\Cc}{{\mathbb C}}

\newcommand{\Zz}{{\mathbb Z}}
\newcommand{\Qq}{{\mathbb Q}}
\newcommand{\Pp}{{\mathbb P}}
\newcommand{\Rr}{{\mathbb R}}
\newcommand{\Nn}{{\mathbb N}}

\newcommand{\OQq}{\overline{{\mathbb Q}}}
\newcommand{\WW}{\mbox{$\mathcal{W}$}}
\newcommand{\HH}{\mbox{$\mathcal{H}$}}

\newcommand{\DD}{\mbox{$\mathcal{D}$}}

\newcommand{\rank}{\operatorname{rank}}

\newcommand{\x}{{\bf x}}
\newcommand{\m}{{\bf m}}

\renewcommand{\leq}{\leqslant}
\renewcommand{\geq}{\geqslant}

\newcommand{\kdots}{,\ldots ,}

\title[Subspace Theorem with polynomials of higher degree]
{A generalization of the Subspace Theorem with polynomials of higher degree}


\author[J.-H.~EVERTSE]{Jan-Hendrik~EVERTSE}
\author[R.~G.~FERRETTI]{Roberto~G.~FERRETTI}
\address{J.-H. Evertse,
Universiteit Leiden, Mathematisch Instituut,
Postbus 9512, 2300 RA Leiden,
The Netherlands}
\email{evertse@math.leidenuniv.nl}
\address{R.G. Ferretti,
Universit\`a della Svizzera Italiana,
Via Buffi 23, CH-6900 Lugano, Switzerland}
\email{roberto.ferretti@lu.unisi.ch}

\begin{document}

\maketitle

\begin{center}
\emph{To Professor Wolfgang Schmidt on his 70th birthday}
\\[0.5cm]
\end{center}

\begin{abstract}
Recently, Corvaja and Zannier \cite[Theorem 3]{CZ} proved an extension
of the Subspace Theorem with polynomials of arbitrary degree
instead of linear forms. Their result states that the set of
solutions in $\Pp^n(K)$ ($K$ number field) of the inequality
being considered
is not Zariski dense.

In this paper we prove, by a different
method, a generalization of their result, in which the solutions
are taken from an arbitrary projective variety $X$ instead of $\Pp^n$.
Further we give a quantitative version, which states
in a precise form that the solutions with
large height
lie in a finite number of proper subvarieties of $X$, with explicit
upper bounds for the number and for the degrees of these subvarieties
(\thmref{th:1.2} below).

We deduce our generalization from a general
result on twisted heights on projective varieties (\thmref{th:4.1}
in \secref{4}).
Our main tools are the quantitative version of the Absolute Parametric
Subspace Theorem by Evertse and Schlickewei \cite[Theorem 1.2]{ES},
as well as a lower bound by Evertse and Ferretti \cite[Theorem 4.1]{EF}
for the normalized Chow weight of a projective variety in terms of its
$m$-th normalized Hilbert weight.
\end{abstract}

\section{Introduction}\label{1}

\begin{prg}
The Subspace Theorem can be stated as follows.
Let $K$ be a number field (assumed to be contained
in some given algebraic closure $\OQq$ of $\Qq$),
$n$ a positive integer, $0<\delta \leq 1$
and $S$ a finite set of places of $K$.
For $v\in S$, let $L_0^{(v)}\kdots L_n^{(v)}$ be linearly
independent linear forms in $\OQq [x_0\kdots x_n]$.
Then the set of solutions $\x\in\Pp^n(K)$ of
\begin{equation}
\label{1.1}
\log \left(
\prod_{v\in S}\prod_{i=0}^n \frac{|L_i^{(v)}(\x )|_v}{\|\x\|_v}
\right)
\leq -(n+1+\delta )h(\x )
\end{equation}
is contained in the union of finitely many proper linear subspaces of $\Pp^n$.

Here,
$h(\cdot )$ denotes the absolute logarithmic height
on $\Pp^n(\OQq )$,
$|\cdot |_v$, $\|\cdot \|_v$ $(v\in S)$ denote
normalized absolute values on $K$ and normalized norms on $K^{n+1}$,
and each
$|\cdot |_v$ has been extended to $\OQq$ (see \prgref{notation} below).
The Subspace Theorem was first proved by Schmidt \cite{S1},\cite{S2}
for the case that $S$ consists
of the archimedean places of $K$, and then later extended by
Schlickewei \cite{Schl} to the general case.
\end{prg}

\begin{prg}
We state a generalization of the Subspace Theorem in which the linear forms
$L_i^{(v)}$ are replaced by homogeneous polynomials of arbitrary degree,
and in which the solutions are taken from an $n$-dimensional
projective subvariety of $\Pp^N$ where $N\geq n\geq 1$.

By a projective subvariety of $\Pp^N$ we mean a geometrically
irreducible Zariski-closed subset
of $\Pp^N$. For a Zariski-closed
subset
$X$ of $\Pp^N$ and for a field $\Omega$, we denote by $X(\Omega )$ the set
of $\Omega$-rational points of $X$.
For homogeneous polynomials
$f_1\kdots f_r$ in the variables $x_0\kdots x_N$
we denote by \\$\{ f_1=0\kdots f_r=0\}$
the Zariski-closed subset of $\Pp^N$ given by $f_1=0\kdots \\f_r=0$.

Then our result reads as follows:
\end{prg}
\vskip 0.3cm
\begin{thm}\label{th:1.1}
Let $K$ be a number field, $S$ a finite set of places of $K$ and $X$ a
projective subvariety of $\Pp^N$ defined over $K$ of dimension $n\geq 1$
and degree $d$.
Let $0<\delta \leq 1$.
Further,
for $v\in S$ let $f_0^{(v)}\kdots f_n^{(v)}$ be a system of homogeneous
polynomials in $\OQq [x_0\kdots x_N]$ such that
\begin{equation}\label{1.3}
X(\OQq )\cap \big\{ f_0^{(v)}=0\kdots f_n^{(v)}=0\big\}\,=\,\emptyset\,\,\,
\mbox{for $v\in S$.}
\end{equation}
Then the set of solutions $\x\in X(K)$ of the inequality
\begin{equation}\label{1.4}
\log\left(
\prod_{v\in S}\prod_{i=0}^{n}
\frac{|f_i^{(v)}(\x )|_v^{1/\deg f_i^{(v)}}}{\| \x \|_v}\right)
\leq -(n+1+\delta )h(\x)
\end{equation}
is contained in a finite union
$\bigcup_{i=1}^u \Big(X\cap \{ G_i =0\}\Big)$,
where $G_1\kdots G_u$ are homogeneous polynomials
in $K[x_0\kdots x_N]$
not vanishing
identically on $X$ of degree at most
\[
(8n+6)(n+2)^2d\Delta^{n+1}\delta^{-1}\quad\text{with }
\Delta :={\rm lcm}\big(\deg f_i^{(v)}:\, v\in S,\, 0\leq i\leq n\big)\, .
\]
\end{thm}
\vskip 0.2cm
It should be noted that if $N=n$, $X=\Pp^n$ and
$f_0^{(v)}\kdots f_n^{(v)}$
are linear forms, then condition \formref{1.3} means precisely that
$f_0^{(v)}\kdots f_n^{(v)}$ are linearly independent.

We give an immediate consequence:
\vskip0.3cm
\begin{cor}\label{co:1.1}
Let $f_0\kdots f_n$ be homogeneous polynomials in $\OQq [x_0\kdots x_n]$
such that
\[
\big\{ \x\in\OQq^{n+1}:\, f_0(\x)=\cdots =f_n(\x )=0\big\}\,=\,\{{\bf 0}\}.
\]
Let $0<\delta \leq 1$. Then the set of solutions
$\x =(x_0\kdots x_n)\in\Zz^{n+1}$ of
\[
\prod_{i=0}^n |f_i(\x )|^{1/deg f_i} \leq \Big(\max_{0\leq i\leq n} |x_i|\Big)^{-\delta}
\]
is contained in some finite union of
hypersurfaces $\{G_1=0\}\cup\cdots\cup \{G_u=0\}$, where each $G_i$
is a homogeneous polynomial in
$\Qq [x_0\kdots x_n]$ of degree at most
$(8n+6)(n+2)^2\Delta^{n+1}\delta^{-1}$ with
$\Delta :={\rm lcm}\big(\deg f_i:\,  0\leq i\leq n\big)$.
\end{cor}
\vskip0.3cm
\begin{prg}
In their paper \cite{FW}, Faltings and W\"{u}stholz introduced
a new method to prove the Subspace Theorem, and gave some examples
showing that their method enables to prove extensions of the
Subspace Theorem with higher degree polynomials instead of linear
forms, and with solutions from an arbitrary projective variety.
Ferretti
\cite{F},\cite{F2} observed the role of Mumford's degree of contact
\cite{M}
(or the Chow weight, see \prgref{Chow} below) in the work of Faltings and W\"{u}stholz
and worked out several other cases. Evertse and Ferretti \cite{EF} showed
that the extensions of the Subspace Theorem as proposed by Faltings
and W\"{u}stholz can be deduced directly from the Subspace Theorem
itself.

Recently,
Corvaja and Zannier \cite[Theorem 3]{CZ} obtained
a result similar to our \thmref{th:1.1} with $X=\Pp^n$.
(More precisely, Corvaja and Zannier gave an essentially equivalent
affine formulation, in which the polynomials $f_i^{(v)}$
need not be homogeneous and in which the solutions ${\bf x}$
have $S$-integer coordinates). In fact, Corvaja and Zannier showed that
the set of solutions of \formref{1.4} is contained in a finite union
of hypersurfaces in $\Pp^n$ and gave some further information
about the structure of these hypersurfaces, on the other hand they did not
provide an explicit bound for their degrees.
Corvaja and Zannier stated their result only for the case
$X=\Pp^n$ but with their
methods this may be extended to the case that $X$
is a complete intersection.
In contrast, our result is valid for arbitrary projective
subvarieties $X$ of $\Pp^N$.

In their paper \cite{CZ}, Corvaja and Zannier proved also
finiteness results for several classes of Diophantine equations.
It is likely, that similar results can be deduced by means
of our approach, but we have not gone into this.
\end{prg}

\begin{prg}\label{notation}
Below we state a quantitative version of \thmref{th:1.1}.
We first introduce the necessary notation.
All number fields considered in this paper are contained
in a given algebraic closure $\OQq$ of $\Qq$.
Let $K$ be a number field and
denote by $G_K$ the Galois group
of $\OQq$ over $K$. For $\x =(x_0\kdots x_N)\in\OQq^{N+1}$, $\sigma\in G_K$
we write $\sigma (\x )= (\sigma (x_0)\kdots \sigma (x_N))$.
Denote by $M_K$
the set of places of $K.$ For $v\in M_K,$ choose an absolute value
$|.|_v$ normalized such that the restriction of $|.|_v$ to $\Qq$
is $|.|^{[K_v:\Rr]/[K:\Qq]}$ if $v$ is archimedean and
$|.|_p^{[K_v:\Qq_p]/[K:\Qq]}$ if $v$ lies above the prime number
$p.$ Here $|.|$ is the ordinary absolute value, and $|.|_p$ is the
$p$-adic absolute value with $|p|_p=p^{-1}.$
These absolute values satisfy the product formula $\prod_{v\in M_K} |x|_v=1$
for $x\in K^*$.

Given $\x =(x_0\kdots x_N)\in K^{N+1}$ we put
$\| \x \|_v:=\max (|x_0|_v\kdots |x_N|_v)$ for $v\in M_K$.
Then
the absolute logarithmic height of $\x$ is defined by
$h(\x )=\log\Big(\prod_{v\in M_K} \| \x \|_v\Big)$.
By the product formula,
$h(\lambda\x )=h(\x)$ for $\lambda\in K^*$. Moreover,
$h(\x)$ depends only on $\x$ and not on the choice of the
particular number field $K$ containing $x_0\kdots x_N$.
Thus, this function
$h$ gives rise to a height on $\Pp^N(\OQq )$.

Given a
system $f_0\kdots f_m$ of polynomials with coefficients in $\OQq$ we define
$h(f_0\kdots f_m):=h({\bf a})$, where ${\bf a}$ is a vector consisting
of the non-zero coefficients of $f_0\kdots f_m$.
Further by $K(f_0\kdots f_m)$ we denote the extension of $K$ generated
by the coefficients of $f_0\kdots f_m$.
The
height of a projective subvariety $X$ of $\Pp^N$ defined over $\OQq$
is defined by $h(X):=h(F_X)$, where $F_X$ is the Chow form
of $X$ (see \prgref{Chow} below).

For every $v\in M_K$
we choose an extension of $|\cdot |_v$ to $\OQq$ (this amounts to extending
$|\cdot |_v$ to the algebraic closure $\overline{K}_v$ of $K_v$ and choosing
an embedding of $\OQq$ into $\overline{K}_v$). Further for $v\in M_K$,
$\x =(x_0\kdots x_N)\in\OQq^{N+1}$ we put
$\|\x \|_v:=\max (|x_0|_v\kdots |x_N|_v)$.
\end{prg}

\begin{prg}
Schmidt \cite{S} was the first to obtain a quantitative version
of the Subspace Theorem, giving an explicit upper bound for the number of
subspaces containing all solutions with `large' height.
Since then his basic result has been improved and generalized
in various directions.
Evertse and Schlickewei \cite[Theorem 3.1]{ES} deduced
a quantitative version of the
Absolute Subspace Theorem,
dealing with solutions in $\Pp^n(\OQq )$ of some absolute extension
of \formref{1.1}. Their result can be stated as follows.
\\[0.2cm]
Let again $K$ be a number field, and $S$
a finite set of places of $K$ of cardinality $s$.
Let $n\geq 1$, $0<\delta\leq 1$. For $v\in S$, let $L_0^{(v)}\kdots L_n^{(v)}$
be linearly independent linear forms in $\OQq [x_0\kdots x_n]$. Put
$\DD :=\prod_{v\in S}|\det (L_0^{(v)}\kdots L_n^{(v)})|_v$ and assume that
$[K(L_i^{(v)}):K]\leq C$ for $v\in S$, $i=0\kdots n$.
Then the set of $\x\in\Pp^n(\OQq )$ with
\[
\log\left(
\DD^{-1}\prod_{v\in S}\prod_{i=0}^n\max_{\sigma\in G_K}
\frac{|L_i^{(v)}(\sigma (\x ))|_v}{\|\sigma (\x )\|_v}\right)
\leq -(n+1+\delta )h(\x )\, ,
\]
\[
h(\x )\geq 9(n+1)\delta^{-1}\log (n+1)+\max\big( h(L_i^{(v)}):\, v\in S,\,
0\leq i\leq n\big)
\]
is contained in the union of not more than
\[
(3n+3)^{(2n+2)s}8^{(n+10)^2}\delta^{-(n+1)s-n-5}\log (4C)\log\log (4C)
\]
proper linear subspaces of $\Pp^n(\OQq )$ which are all defined over $K$.
\\[0.2cm]
Typically, the lower bound for $h(\x )$ depends on the linear
forms $L_i^{(v)}$, while the upper bound for the number
of subspaces does not depend on the $L_i^{(v)}$.
\end{prg}

\begin{prg}
We now state an analogue for inequalities with higher degree polynomials
instead of linear forms. We first list some notation:
\\[0.2cm]
$\delta$ is a real with $0<\delta \leq 1$,
$K$ is a number field, $S$ is a finite set of places of $K$ of cardinality
$s$, $X$
is a projective subvariety of $\Pp^N$ defined over $K$
of dimension $n\geq 1$ and degree $d$,
$f_0^{(v)}\kdots f_n^{(v)}$ $(v\in S)$
are systems of homogeneous polynomials in
$\OQq [x_0\kdots x_N]$,
\begin{eqnarray}
\label{1.5a}
&&
\left\{
\begin{array}{l}
C:=\max\big( [K(f_i^{(v)}):K]:\, v\in S,\, i=0\kdots n\big),
\\
\Delta :={\rm lcm}\big(\deg f_i^{(v)}:\, v\in S,\, i=0\kdots n\big),
\end{array}\right.
\\[0.4cm]
\label{1.5b}
&&
\left\{
\begin{array}{l}
A_1:= (20n\delta^{-1})^{(n+1)s}\cdot
\exp\Big( 2^{12n+16}n^{4n}\delta^{-2n}d^{2n+2}\Delta^{n(2n+2)}\Big)\cdot
\\
\qquad\qquad\qquad\qquad
\qquad\qquad\qquad\qquad\quad
\cdot\log (4C)\log\log (4C),
\\[0.2cm]
A_2:=(8n+6)(n+2)^2d\Delta^{n+1}\delta^{-1},
\\[0.2cm]
A_3:= \exp\Big( 2^{6n+20}n^{2n+3}\delta^{-n-1}d^{n+2}\Delta^{n(n+2)}\log (2Cs)\Big),
\\[0.2cm]
H:= \log (2N) +h(X)+
\max\big( h(1,f_i^{(v)}):\, v\in S,\, 0\leq i\leq n\big).
\end{array}\right.
\end{eqnarray}
\end{prg}
\vskip 0.2cm
\begin{thm}\label{th:1.2}
Assume that
\\[0.2cm]
{\rm (\ref{1.3})}$\qquad\qquad\,$
$X(\OQq )\cap \big\{ f_0^{(v)}=0\kdots f_n^{(v)}=0\big\}\,=\,\emptyset\quad$
for $v\in S$.
\\[0.2cm]
Then there are homogeneous polynomials
$G_1\kdots G_u\in K[x_0\kdots x_N]$ with
\[
u\leq A_1\, ,\quad \deg G_i\leq A_2\,\,\,\mbox{for $i=1\kdots u$}
\]
which do not vanish identically on $X$, such that
the set of $\x\in X(\OQq )$ with
\begin{equation}\label{1.6}
\log\left(
\prod_{v\in S}\prod_{i=0}^{n}
\max_{\sigma\in G_K}
\frac{|f_i^{(v)}(\sigma (\x ))|_v^{1/\deg f_i^{(v)}}}{\|\sigma (\x )\|_v}\right)
\leq -(n+1+\delta )h(\x)\, ,
\end{equation}
\begin{equation}\label{1.7}
h(\x )\geq A_3\cdot H
\end{equation}
is contained in $\bigcup_{i=1}^u\Big( X\cap \{G_i=0\}\Big)$.
\end{thm}
\vskip 0.3cm
Clearly, the bounds in \thmref{th:1.2} are much worse than those in the
result of Evertse and Schlickewei.
It would be very interesting if one could replace $A_1,A_3$ by quantities
which are at most exponential in (some power of) $n$ and which are polynomial
in $\delta^{-1},d,\Delta$. Further, we do not know
whether the dependence of $A_2$ on $\delta$ is needed.

\begin{prg}
Our starting point is a result for twisted heights
on $\Pp^n$ (a quantitative version of the
Absolute Parametric
Subspace Theorem), due to Evertse and Schlickewei
\cite[Theorem 2.1]{ES} (see also \propref{pr:5.1} in
\secref{5} below).
From this, we deduce
an analogous result for
twisted heights on arbitrary projective varieties;
the statement of this result is in \secref{4} (\thmref{th:4.1})
and its proof in \secref{5}.
The proof involves some arguments
from Evertse and Ferretti \cite{EF}, in particular
an explicit lower bound of the normalized Chow weight of a projective
variety in terms of the $m$-th normalized Hilbert weight
of that variety.
In \secref{2} we give some height estimates;
here we use heavily R\'{e}mond's expos\'{e} \cite{R2}.
Then in \secref{3} we deduce \thmref{th:1.2}.
Using that $\Pp^N(K)$ has only finitely many points with height below
any given bound,
\thmref{th:1.1}
follows at once from \thmref{th:1.2}.
\end{prg}
\vskip0.5cm
\section{Twisted heights}\label{4}

\begin{prg}
The quantitative version of the Absolute Parametric
Subspace Theorem of Evertse and Schlickewei mentioned in the previous
section
deals with a class of twisted heights defined
on $\Pp^n(\OQq )$ parametrized
by a real $Q\geq 1$. Roughly speaking, this result states that there
are a finite number of proper linear subspaces of $\Pp^n$ such that for every
sufficiently large $Q$, the set of points in $\Pp^n(\OQq )$ with small
$Q$-height is contained in one of these subspaces.
\thmref{th:4.1} stated below is an analogue
in which the points are taken from an arbitrary projective variety
instead of $\Pp^n$.
Loosely speaking, \thmref{th:1.2} is proved
by defining a suitable finite morphism $\varphi$ from
$X$ to a projective variety $Y\subset\Pp^R$ and
a finite number of classes of twisted heights on $Y$ as above,
and applying \thmref{th:4.1} to each of these classes.
\end{prg}

\begin{prg}
Let $K$ be a number field. For finite extensions
of $K$ we define normalized absolute values similarly
as for $K$. Thus, if $L$ is a finite extension of $K$,
$w$ is a place of $L$, and $v$ is the place of $K$ lying
below $w$, then
\begin{equation}\label{4.1}
|x|_w=|x|_v^{d(w|v)}\,\,\,\mbox{for $x\in K$, with }
d(w|v):=\frac{[L_w:K_v]}{[L:K]},
\end{equation}
where $K_v$, $L_w$ denote the completions at $v,w$,
respectively.

We denote points on $\Pp^R$ by ${\bf y}=(y_0\kdots y_R)$.
For $v\in M_K$, let ${\bf c}_v=(c_{0v}\kdots c_{Rv})$ be a
tuple of reals such that $c_{0v}=\cdots =c_{Rv}=0$ for all
but finitely many places $v\in M_K$ and put
${\bf c}=({\bf c}_v:\, v\in M_K)$.
Further, let $Q$ be a real $\geq 1$. We define a twisted height
on $\Pp^R(\OQq )$ as follows. First put
\[
H_{Q,{\bf c}}({\bf y}):=
\prod_{v\in M_K}\max_{0\leq i\leq R}\Big( |y_i|_vQ^{c_{iv}}\Big)
\quad\mbox{for ${\bf y}=(y_0\kdots y_R)\in\Pp^R(K)$;}
\]
by the product formula, this is well-defined on $\Pp^R(K)$.
For any finite extension $L$ of $K$ we put
\begin{equation}\label{4.2}
c_{iw}:=c_{iv}\cdot d(w|v)\quad\mbox{for $w\in M_L$,}
\end{equation}
where $M_L$ is the set of places of $L$ and $v$ the place of
$K$ lying below $w$. Then for ${\bf y}\in\Pp^R(\OQq )$, we define
\begin{equation}\label{4.3}
H_{Q,{\bf c}}({\bf y}):=
\prod_{w\in M_L}\max_{0\leq i\leq R}\Big( |y_i|_wQ^{c_{iw}}\Big)
\end{equation}
where $L$ is any finite extension of $K$ such that
${\bf y}\in\Pp^R(L)$.
In view of \formref{4.1} this definition does not depend on $L$.
\end{prg}

\begin{prg}\label{Chow}
Let $Y$ be a (by definition irreducible) projective subvariety of $\Pp^R$ of
dimension $n$ and degree $D$, defined over $K$. We recall that
there is an
up to a constant factor unique
polynomial $F_Y({\bf u}^{(0)}\kdots {\bf u}^{(n)})$ with
coefficients in $K$ in blocks of variables ${\bf
u}^{(0)}=(u_0^{(0)}\kdots u_R^{(0)})$$\kdots$ ${\bf
u}^{(n)}=(u_0^{(n)}\kdots u_R^{(n)})$, called the \emph{Chow form} of $Y$,
with the following
properties:
\\[0.2cm]
$F_Y$ is irreducible over $\OQq$; $F_Y$ is homogeneous
in each block ${\bf u}^{(h)}$ ($h=0\kdots n$); and $F_Y({\bf
u}^{(0)}\kdots {\bf u}^{(n)})=0$ if and only if $Y$ and the
hyperplanes $\sum_{i=0}^R u_i^{(h)}y_i=0$ $(h=0\kdots n)$ have a
$\OQq$-rational point in common.
\\[0.2cm]
It is well-known that the degree
of $F_Y$ in each block ${\bf u}^{(h)}$ is $D$.

Let ${\bf c}=(c_0\kdots c_R)$ be a tuple of reals.
Introduce an auxiliary variable $t$ and substitute $t^{c_i}u_i^{(h)}$
for $u_i^{(h)}$ in $F_Y$ for $h=0\kdots n$, $i=0\kdots R$.
Thus we obtain an expression
\begin{eqnarray}
\label{4.4}
&&
F_Y(t^{c_0}u_0^{(0)}\kdots t^{c_R}u_R^{(0)};\ldots;t^{c_0}u_0^{(n)}\kdots t^{c_R}u_R^{(n)})\\
\nonumber
&&\qquad\qquad =  t^{e_0}G_0({\bf u}^{(0)}\kdots {\bf u}^{(n)})+\cdots +
t^{e_r}G_r({\bf u}^{(0)}\kdots {\bf u}^{(n)}),
\end{eqnarray}
with $G_0\kdots G_r\in K [{\bf u}^{(0)}\kdots {\bf u}^{(n)}]$ and
$e_0>e_1>\cdots>e_r$. We now define the {\it Chow weight} of $Y$
with respect to ${\bf c}$
\footnote{The Chow weight is closely related to the \emph{degree of contact}
earlier introduced by Mumford \cite{M}. Roughly speaking, the degree of contact
of $Y$ with respect to ${\bf c}$ is defined for integer tuples ${\bf c}$
and it is equal to $e_r$ instead of $e_0$.}
by
\begin{equation}\label{4.5}
e_Y({\bf c}):=e_0.
\end{equation}
\end{prg}

\begin{prg}
We formulate our main result for twisted heights.
Below, $Y$ is a projective subvariety of $\Pp^R$ of dimension
$n\geq 1$ and degree $D$, defined over $K$, and
${\bf c}_v=(c_{0v}\kdots c_{Rv})$ $(v\in M_K)$ are tuples
of reals such that
\begin{eqnarray}
\label{4.6}
&&c_{iv}\geq 0\,\,\,\mbox{for $v\in M_K$, $i=0\kdots R$;}
\\
\label{4.6b}
&&c_{0v}=\cdots =c_{Rv}=0\,\,\,\mbox{for all but finitely many $v\in M_K$;}
\\
\label{4.6c}
&&\sum_{v\in M_K}\max (c_{0v}\kdots c_{Rv})\leq 1.
\end{eqnarray}
Put
\begin{equation}\label{4.6a}
E_Y({\bf c}) :=\frac{1}{(n+1)D}\left(\sum_{v\in M_K} e_Y({\bf c}_v)\right)\, .
\end{equation}
Further, let $0<\delta \leq 1$, and put
\begin{equation}\label{4.7}
\left\{ \quad
\begin{array}{l}
B_1:=\exp\Big( 2^{10n+4}\delta^{-2n}D^{2n+2}\Big)
\cdot\log (4R)\log\log (4R),
\\
B_2:= (4n+3)D\delta^{-1},
\\
B_3:= \exp\Big( 2^{5n+4}\delta^{-n-1}D^{n+2}\log (4R)\Big).
\end{array}\right.
\end{equation}
\vskip0.2cm
\begin{thm}\label{th:4.1}
There are homogeneous polynomials
$F_1\kdots F_t\in\\ K[y_0\kdots y_R]$ with
\[
t\leq B_1,\quad \deg F_i\leq B_2\,\,\,\mbox{for $i=1\kdots t$,}
\]
which do not vanish identically on $Y$, such that for every real
number $Q$ with
\[
\log Q\geq B_3\cdot (h(Y)+1)
\]
there is $F_i\in\{ F_1\kdots F_t\}$ with
\begin{equation}\label{4.8}
\left\{ {\bf y}\in Y(\OQq ):\, H_{Q,{\bf c}}({\bf y})\leq
Q^{E_Y({\bf c})-\delta}\right\} \subset Y\cap\big\{F_i=0\big\}\, .
\end{equation}
\end{thm}
\end{prg}
\vskip0.5cm
\section{Proof of \thmref{th:4.1}}\label{5}

\begin{prg}
We first recall the quantitative version of the Absolute Parametric Subspace
Theorem of Evertse and Schlickewei.
As before, $K$ is an algebraic number field and $R,n$ are integers
with $R\geq n\geq 1$.
We denote the coordinates on $\Pp^n$ by
$(x_0\kdots x_n)$. Given an index set $I=\{ i_0\kdots i_n\}$
with $i_0<\cdots <i_n$
and linear forms $L_j=\sum_{i=0}^n a_{ij}x_i$ $(j\in I)$
we write $\det (L_j:\, j\in I):= \det (a_{i,i_j})_{i,j=0\kdots n}$.

Let $L_0\kdots L_R$ be linear forms in $K[x_0\kdots x_n]$ with
$\rank\{ L_0\kdots L_R\}=n+1$. Further, let
$I_v$ $(v\in M_K)$ be subsets of
$\{ 0\kdots R\}$ of cardinality $n+1$ such that
\begin{equation}\label{5.1}
\rank\{ L_i:\, i\in I_v\}=n+1\quad\mbox{for $v\in M_K$.}
\end{equation}
Define
\begin{equation}\label{5.3}
\HH := \prod_{v\in M_K}\max_I |\det (L_i:\, i\in I)|_v\, ,\quad
\DD := \prod_{v\in M_K} |\det (L_i:\, i\in I_v)|_v\, ;
\end{equation}
here the maximum is taken over all subsets $I$ of
$\{ 0\kdots R\}$ of cardinality $n+1$. According to \cite[Lemma 7.2]{EF}
we have
\begin{equation}\label{5.3a}
\DD\geq \HH ^{1-\binom{R+1}{n+1}}\, .
\end{equation}
Let ${\bf d}_v=(d_{iv}:\, i\in I_v)$ $(v\in M_K)$
be tuples of reals such that
\begin{eqnarray}
\label{5.2}
&&d_{iv}=0\,\,\,\mbox{for $i\in I_v$ and for all but finitely many $v\in M_K$,}
\\[0.2cm]
\label{5.2a}
&&\sum_{v\in M_K}\sum_{i\in I_v} d_{iv}=0,
\\[0.2cm]
\label{5.2b}
&&\sum_{v\in M_K}\max (d_{iv}:\, i\in I_v)\leq 1
\end{eqnarray}
and write ${\bf d}=({\bf d}_v:\, v\in M_K)$.
\\[0.2cm]

We define a twisted height on $\Pp^n(\OQq )$ as follows. For any
real number $Q\geq 1$ we first put
\[
H^*_{Q,{\bf d}}({\bf x})=
\prod_{v\in M_K}\left(\max_{i\in I_v} |L_i({\bf x})|_vQ^{-d_{iv}}\right)
\quad\mbox{for ${\bf x}\in\Pp^n(K)$.}
\]
More generally, if $L$ is any finite extension of $K$, put
\begin{equation}\label{5.6}
d_{iw}:= d(w|v)d_{iv},\,\,\, I_w:= I_v
\end{equation}
where $v$ is the place of $K$ lying below $w$.
Then for ${\bf x}\in\Pp^n(\OQq )$ we define
\begin{equation}\label{5.7}
H^*_{Q,{\bf d}}({\bf x})=
\prod_{w\in M_L}\left(\max_{i\in I_w} |L_i({\bf x})|_vQ^{-d_{iw}}\right)
\end{equation}
where $L$ is any finite extension of $K$ such that ${\bf x}\in\Pp^n(L)$.
This is independent of the choice of $L$.

Now the result of Evertse and Schlickewei \cite[Theorem 2.1]{ES}
is as follows:
\end{prg}
\vskip 0.3cm
\begin{prop}\label{pr:5.1}
Let $I_v$ $(v\in M_K)$, ${\bf d}=({\bf d}_v:\, v\in M_K)$,
satisfy \formref{5.1}, \formref{5.2}, respectively,
and let $0<\varepsilon \leq 1$.
\\
There are proper linear subspaces $T_1\kdots T_t$
of $\Pp^n$, defined over $K$, with
\begin{equation}\label{5.8}
t\leq 4^{(n+9)^2}\varepsilon^{-n-5}\log (3R)\log\log (3R),
\end{equation}
such that for every real number $Q$ with
\begin{equation}\label{5.9}
Q\geq \max\left( \HH^{1/\binom{R+1}{n+1}}, (n+1)^{2/\varepsilon}\right)
\end{equation}
there is $T_i\in\{ T_1\kdots T_t\}$ with
\begin{equation}\label{5.10}
\{ {\bf x}\in\Pp^n(\OQq ):\,
H^*_{Q,{\bf d}}({\bf x})\leq
\DD^{1/(n+1)}Q^{-\varepsilon}\}\subset T_i\, .
\end{equation}
\end{prop}
\vskip0.4cm

\begin{prg}
We recall some results from \cite{EF}. As in \secref{4}, we denote
the coordinates on $\Pp^R$ by $(y_0\kdots y_R)$. Let $Y$ be a
projective variety of $\Pp^R$ defined over $K$ of dimension $n$
and degree $D$. Let $I_Y$ be the prime ideal of $Y$, i.e. the
ideal of polynomials from $\OQq [y_0\kdots y_R]$ vanishing
identically on $Y$. For $m\in\Nn$, denote by $\OQq [y_0\kdots
y_R]_m$ the vector space of homogeneous polynomials in $\OQq
[y_0\kdots y_R]$ of degree $m$, and put $(I_Y)_m:=$$\OQq
[y_0\kdots y_R]_m\cap I_Y$. Then the Hilbert function of $Y$ is
defined by
\[
H_Y(m):=\dim_{\OQq}\Big(\OQq [y_0\kdots y_R]_m/(I_Y)_m\Big)\, .
\]
The scalar product of ${\bf a}=(a_0\kdots a_R),\,\,{\bf b}=(b_0\kdots b_R)
\in\Rr^{R+1}$ is given by ${\bf a}\cdot {\bf b}:=a_0b_0+\cdots +a_Rb_R$.
For ${\bf a}=(a_0\kdots a_R)\in(\Zz_{\geq 0})^{R+1}$,
denote by ${\bf y}^{{\bf a}}$ the monomial
$y_0^{a_0}\cdots y_R^{a_R}$.
Then the $m$-th Hilbert weight of $Y$ with respect to a tuple
${\bf c}=(c_0\kdots c_R)\in\Rr^{R+1}$ is defined by
\begin{equation}\label{5.14}
s_Y(m,{\bf c})
:=\max\left( \sum_{i=1}^{H_Y(m)} {\bf a}_i\cdot {\bf c}\right)\, ,
\end{equation}
where the maximum is taken over all sets of monomials
$\{ {\bf y}^{{\bf a}_1}\kdots {\bf y}^{{\bf a}_{H_Y(m)}}\}$,
whose residue classes
modulo $(I_Y)_m$ form a basis of $\OQq [y_0\kdots y_R]_m/(I_Y)_m$.
\end{prg}
We recall Evertse and Ferretti \cite[Theorem 4.1]{EF}:
\vskip0.3cm
\begin{prop}\label{pr:5.2}
Let ${\bf c}=(c_0\kdots c_R)$ be a tuple of non-negative reals. Let $m>D$
be an integer. Then
\begin{equation}\label{5.15}
\mbox{$\frac{1}{mH_Y(m)}\cdot s_Y(m,{\bf c})\geq
\frac{1}{(n+1)D}\cdot e_Y({\bf c})-\frac{(2n+1)D}{m}
\cdot \max (c_0\kdots c_R)\, .$}
\end{equation}
\end{prop}
\vskip0.4cm

Let $m$ be a positive integer. Put
\[
n_m:=H_Y(m)-1,\quad R_m:=\mbox{$\binom{R+m}{m}$}-1,
\]
and let ${\bf y}^{{\bf a}_0}\kdots {\bf y}^{{\bf a}_{R_m}}$ be the monomials
of degree $m$ in $y_0\kdots y_R$, in some order. Denote by $\varphi_m$
the Veronese map of degree $m$,
${\bf y}\mapsto ({\bf y}^{{\bf a}_0}\kdots {\bf y}^{{\bf a}_{R_m}})$.
Lastly, denote by $Y_m$ the smallest linear subspace of $\Pp^{R_m}$
containing $\varphi_m(Y)$.
\vskip0.3cm
\begin{lemma}\label{le:5.3}
{\rm (i)} $Y_m$ is defined over $K$;\\
{\rm (ii)} $\dim Y_m =n_m\leq D\binom{m+n}{n}$;\\
{\rm (iii)}
$h(Y_m)\leq Dm\binom{m+n}{n}\Big( D^{-1}h(Y)+(3n+4)\log (R+1)\Big)$.
\end{lemma}
\vskip0.2cm
\emph{Proof.} (i),(iii) \cite[Lemma 8.3]{EF};
(ii) Chardin \cite[Th\'{e}or\`{e}me 1]{Ch}.\qed
\vskip0.5cm
\begin{prg}
Let ${\bf c}_v\in\Rr^R$ $(v\in M_K)$ be tuples with \formref{4.6}.
For a suitable value of $m$, we link the twisted height $H_{Q,{\bf c}}$
from \thmref{th:4.1} to a twisted height
on $\Pp^{n_m}$ to which \propref{pr:5.1} is applicable.
Put
\begin{equation}\label{5.16}
m:=[(4n+3)D\delta^{-1}]\, .
\end{equation}
Then by \propref{pr:5.2} and \formref{4.6} we have
\begin{equation}\label{5.17}
\frac{1}{mH_Y(m)}\cdot\left(\sum_{v\in M_K} s_Y(m,{\bf c}_v)\right)\geq
\frac{1}{(n+1)D}\cdot
\left(\sum_{v\in M_K} e_Y({\bf c}_v)\right)-\frac{\delta}{2}\, .
\end{equation}

Denote as before the coordinates on $\Pp^R$ by ${\bf y}=(y_0\kdots y_R)$,
those on $\Pp^{n_m}=\Pp^{H_Y(m)-1}$ by ${\bf x}=(x_0\kdots x_{n_m})$,
and those on $\Pp^{R_m}=\Pp^{\binom{R+m}{m}-1}$ by
${\bf z}=(z_0\kdots z_{R_m})$.
Since $Y_m$ is an $n_m$-dimensional linear subspace of $\Pp^{R_m}$ defined
over $K$, there are linear forms $L_0\kdots L_{R_m}\in K[x_0\kdots x_{n_m}]$
such that the map
\[
\psi_m :\, {\bf x}
\mapsto (L_0({\bf x})\kdots L_{R_m}({\bf x}))
\]
is a linear isomorphism from $\Pp^{n_m}$ to $Y_m$.
Thus, $\psi_m^{-1}\varphi_m$
is an injective map from $Y$ into $\Pp^{n_m}$.

For $v\in M_K$ there is a subset $I_v$ of $\{ 0\kdots R_m\}$ of cardinality
$n_m+1=H_Y(m)$ such that $\{ {\bf y}^{{\bf a}_i}:\, i\in I_v\}$ is a basis
of $\OQq [y_0\kdots y_R]_m/(I_Y)_m$ and
\begin{equation}\label{5.x}
s_Y(m,{\bf c}_v)=\sum_{i\in I_v} {\bf a}_i\cdot {\bf c}_v\, .
\end{equation}
Now define the tuples ${\bf d}_v=(d_{iv}\, ,i\in I_v)$ $(v\in M_K)$ by
\begin{eqnarray}\label{5.20}
d_{iv}&=&-\frac{1}{m}\cdot{\bf a}_i\cdot{\bf c}_v+
\frac{1}{m(n_m+1)}\left(\sum_{j\in I_v} {\bf a}_j\cdot {\bf c}_v\right)
\\
\nonumber &=& -\frac{1}{m}\cdot{\bf a}_i\cdot{\bf c}_v+
\frac{1}{mH_Y(m)}\cdot s_Y(m,{\bf c}_v) \, ,
\end{eqnarray}
and put ${\bf d}=({\bf d}_v:\, v\in M_K)$.
Similarly to \formref{5.3}
we define
\[
\HH :=\prod_{v\in M_K}\max_I |\det (L_i:\, i\in I)|_v,\quad
\DD :=\prod_{v\in M_K} |\det (L_i:\, i\in I_v)|_v\, ,
\]
where the maximum is taken over all subsets $I$ of $\{ 0\kdots R_m\}$
of cardinality $n_m+1$.
Then by, e.g., \cite[page 1300]{EF} we have
\begin{equation}\label{5.21}
\log\HH =h(Y_m)\, .
\end{equation}

We define in a usual manner a twisted height
on $\Pp^{n_m}(\OQq )$ by putting
\[
H^*_{Q,{\bf d}}({\bf x})=\prod_{w\in M_L}\max_{i\in I_w} \big(|L_i({\bf x})|_wQ^{-d_{iw}}\big)
\]
for ${\bf x}\in\Pp^{n_m}(\OQq )$, where $L$ is any finite
extension of $K$ such that ${\bf x}\in\Pp^{n_m}(L)$, $Q\geq 1$ is
a real number, and $d_{iw}=d(w|v)d_{iv}$, $I_w=I_v$  with $v$ the
place of $K$ below $w$. It follows at once from \formref{4.6b}
that $d_{iv}=0$ for all but finitely many $v$ and for $i\in I_v$.
Therefore this height is well-defined.
\end{prg}
\vskip0.3cm
\begin{lemma}\label{le:5.4}
Assume that
\begin{equation}\label{5.21c}
Q\geq \DD^{6/\delta m(n_m+1)}\, .
\end{equation}
Let ${\bf y}\in Y(\OQq )$ be such that
\begin{equation}\label{5.21b}
H_{Q,{\bf c}}({\bf y})
\leq
Q^{E_Y({\bf c})-\delta},
\end{equation}
where
$E_Y({\bf c})=\frac{1}{(n+1)D}\left(\sum_{v\in M_K} e_Y({\bf c}_v)\right)$.
Let ${\bf x}=\psi_m^{-1}\varphi_m({\bf y})$.
Then
\begin{equation}\label{5.21x}
H^*_{Q^m,{\bf d}}({\bf x})\leq \DD^{1/(n_m+1)}(Q^m)^{-\delta /3}\, .
\end{equation}
\end{lemma}
\vskip0.2cm
\emph{Proof.} Put
$s_v:=\frac{1}{mH_Y(m)}s_Y(m,{\bf c}_v)$, $s:=\sum_{v\in M_K} s_v$.
We first show that
\begin{equation}\label{5.22}
H^*_{Q^m,{\bf d}}({\bf x})\leq Q^{-ms}\big(H_{Q,{\bf c}}({\bf y})\big)^m .
\end{equation}
Take a finite extension $L$ of $K$ such that ${\bf y}\in Y(L)$.
We have ${\bf x}\in\Pp^{n_m}(L)$ and
$L_i({\bf x})={\bf y}^{{\bf a}_i}$ for $i=0\kdots R_m$.
So for $w\in M_L$ we have (putting
$s_w:=d(w|v)s_v$, with $v$ the place of $K$ below $w$),
\begin{eqnarray*}
&&\max_{i\in I_w} \big( |L_i({\bf x})|_w (Q^m)^{-d_{iw}}\big)
=\max_{i\in I_w}
\big(|{\bf y}^{{\bf a}_i}|_w Q^{{\bf a}_i\cdot {\bf c}_w-ms_w}\big)
\\
&&\quad\leq \max_{i=0\kdots R_m}
\big( |{\bf y}^{{\bf a}_i} |_w Q^{{\bf a}_i\cdot {\bf c}_w-ms_w}\big)
\leq \left( Q^{-s_w}\max_{i=0\kdots R} \big(|y_i|_wQ^{c_{iw}}\big)\right)^m\, .
\end{eqnarray*}
By taking the product over all $w\in M_L$, \formref{5.22} follows.

Now a successive application of \formref{5.21c},
\formref{5.22}, \formref{5.21b}, \formref{5.17} gives
\[
H^*_{Q^m,{\bf d}}({\bf x}) \leq \DD^{1/(n_m+1)}Q^{m\delta /6}\cdot
Q^{-ms}Q^{mE_Y({\bf c})-m\delta}
\leq \DD^{1/(n_m+1)}(Q^m)^{-\delta/3}\, .
\]
\qed
\vskip0.4cm

\begin{prg}
To complete the proof of \thmref{th:4.1}
we apply \propref{pr:5.1} to \formref{5.21x}; that is,
we apply \propref{pr:5.1} 
with $n=n_m$, $R=R_m$, $\varepsilon=\delta/3$,
and with $Q^m$ in place of $Q$. For the moment we assume
\begin{equation}\label{5.21a}
\log Q\geq \frac{6}{(n_m+1)m\delta}(R_m+1)^{n_m+1}(h(Y_m)+1)\, .
\end{equation}
Notice that this is precisely
\formref{5.9} with
$R=R_m,n=n_m,\varepsilon=\delta/3$ and with $Q^m$ in place of $Q$.

We have to verify that
\formref{5.1}, \formref{5.2}, \formref{5.2a}, \formref{5.2b}
are satisfied with $n_m,R_m$ in place of $n,R$. First,
\formref{5.1} follows at once
from the definition of $I_v$ and the fact that $\psi_m$ is a linear isomorphism.
Second, \formref{5.2} follows from \formref{4.6b}.
Third, \formref{5.2a} follows from \formref{5.20}, \formref{5.x}.
Last, \formref{5.2b} is  consequence of \formref{4.6}, \formref{4.6c}
and the fact that $\frac{1}{mH_Y(m)}\cdot s_Y(m,{\bf c}_v)$ can be expressed
as a maximum of linear forms in $c_{0v}\kdots c_{Rv}$, whose
coefficients are non-negative and have sum equal to $1$.

Thus, there are proper linear subspaces $T_1\kdots T_t$ of
$\Pp^{n_m}$, defined over $K$, with
\begin{equation}\label{5.23}
t\leq 4^{(n_m+9)^2}(3/\delta )^{n_m+5}\log (3R_m)\log\log (3R_m)
\end{equation}
such that for every $Q$ with \formref{5.21a}
there is $T_i\in\{ T_1\kdots T_t\}$ with
\[
\{ {\bf x}\in\Pp^{n_m}(\OQq ):\, H^*_{Q^m,{\bf d}}({\bf x})\leq
\DD^{1/(n_m+1)}(Q^m)^{-\delta /3}\}\subset T_i\, .
\]

For each space $T_i$ there is a linear form
$L_i\in K[z_0\kdots z_{R_m}]$ vanishing identically on $\psi_m(T_i)$
but not on $Y_m$. Since by definition, $Y_m$ is the smallest linear
subvariety of $\Pp^{R_m}$ containing $\varphi_m(Y)$, the linear form
$L_i$ does not
vanish identically on $\varphi_m(Y)$. Replacing in $L_i$ the coordinate
$z_j$ by ${\bf y}^{{\bf a}_j}$ for $j=0\kdots R_m$,
we obtain a homogeneous polynomial
$F_i\in K[y_0\kdots y_R]$ of degree $m$ not vanishing identically on
$Y$ such that if ${\bf x}=\psi_m^{-1}\varphi_m({\bf y})\in T_i$, then
$F_i({\bf y})=0$.

It is easily seen that assumption \formref{5.21a}, together with
\formref{5.21} and \formref{5.3a}, implies \formref{5.21c};
hence \lemref{le:5.4} is applicable.
Thus,
we infer that there are homogeneous polynomials
$F_1\kdots F_t\in K[y_0\kdots y_R]$
of degree $m$,
with $t$ satisfying \formref{5.23}, such that for every $Q$
with \formref{5.21a}
there is $F_i\in \{ F_1\kdots F_t\}$ with
\[
\{ {\bf y}\in Y(\OQq ):\, H_{Q,{\bf c}}({\bf y})\leq
Q^{E_Y({\bf c})-\delta}\}
\subset Y\cap\big\{F_i=0\big\}\, .
\]
By \formref{5.16} we have $m\leq (4n+3)D\delta^{-1}$,
which is the quantity $B_2$ from \formref{4.7}.
So to complete the proof of \thmref{4.1},
it suffices to show that the right-hand side of \formref{5.23} is at most
$B_1$ and that the right-hand side of \formref{5.21a} is at most
$B_3\cdot (h(Y)+1)$, where $B_1,B_3$ are given by \formref{4.7}.

Using $m\geq 7$ and the inequality
\begin{equation}\label{5.binom}
\binom{x+y}{y}\leq\frac{(x+y)^{x+y}}{x^xy^y}=\Big( 1+\frac{y}{x}\Big)^x
\cdot \Big( 1+\frac{x}{y}\Big)^y
\leq \Big(e\big(1+\frac{x}{y}\big)\Big)^y
\end{equation}
for positive integers $x,y$, we infer
\begin{equation}
\label{5.28}
R_m =\binom{R+m}{m}-1
\leq \Big( e\big( 1+\frac{R}{m}\big)\Big)^m\leq (4R)^m\, .
\end{equation}
So by \formref{5.16},
\begin{eqnarray*}
\log (3R_m)\log\log (3R_m) &\leq& 2m^2\log (4R)\log\log (4R)
\\
&\leq& 2(8n+6)^2D^2\delta^{-2}\log (4R)\log\log (4R)\, .
\end{eqnarray*}
Further, by \lemref{le:5.3}, (ii),
\begin{eqnarray}\label{5.29}
n_m &\leq& D\binom{m+n}{n}\leq D\big( e(1+\frac{m}{n})\big)^n
\\
\nonumber
&\leq& D\big( e(1+7D\delta^{-1})\big)^n\leq 2^{5n}\delta^{-n}D^{n+1}\, .
\end{eqnarray}
Hence the right-hand side of \formref{5.23} is at most
\begin{eqnarray*}
&&4^{(2^{5n}\delta^{-n}D^{n+1}+9)^2}(3\delta^{-1})^{2^{5n}\delta^{-n}D^{n+1}+5}
\cdot
\\
&&\qquad\qquad\qquad
\cdot 2(8n+6)^2D^2\delta^{-2}\log (4R)\log\log (4R)
\\
&\leq&
\exp\Big(2^{10n+4}\delta^{-2n}D^{2n+2}\Big)\cdot\log (4R)\log\log (4R)=B_1\, ,
\end{eqnarray*}
while by \lemref{le:5.3}, \formref{5.16}, \formref{5.28}, \formref{5.29},
the right-hand side of \formref{5.21a} is at most
\begin{eqnarray*}
&&\frac{6}{(n_m+1)m\delta}\Big( (4R)^m+1\Big)^{n_m+1}\cdot
\\
&&\qquad\qquad\qquad
\cdot \Big( 1+Dm\binom{m+n}{n}\big( D^{-1}h(Y)+(3n+4)\log (R+1)\Big)
\\
&\leq&
\delta^{-1}\Big( (4R)^{(4n+3)D\delta^{-1}}+1\Big)^{2^{5n}\delta^{-n}D^{n+1}+1}
\cdot
\\
&&\qquad\qquad\qquad
\cdot
2^{5n}\delta^{-n}D^{n+1}(3n+1)\log(R+1)\cdot (h(Y)+1)
\\
&<&
\exp\Big(
2^{5n+4}\delta^{-n-1}D^{n+2}\log (4R)\Big)\cdot (h(Y)+1)=B_3\cdot (h(Y)+1)\, .
\end{eqnarray*}
This completes the proof of \thmref{th:4.1}.\qed
\end{prg}

\vskip0.5cm
\section{Height estimates}\label{2}

\begin{prg}
In this section we compute some height estimates, using
R\'{e}mond's paper \cite{R2}.

Let $K$ be a
number field. Denote as before the set of places of $K$ by $M_K$,
and denote the sets of archimedean and non-archimedean places of $K$
by $M_K^{\infty}$ and $M_K^0$, respectively. We use the normalized
absolute values $|\cdot |_v$ introduced in \prgref{notation}.
Recall that for each of these absolute values we have chosen an
extension to $\OQq$.
In particular, for each $v\in M_K^{\infty}$
there is an isomorphic embedding $\sigma_v:\, \OQq\hookrightarrow\Cc$ such that
$|x|_v=|\sigma_v (x)|^{[K_v:\Rr ]/[K:\Qq  ]}$ for $x\in\OQq$.

We represent
polynomials as
$f=\sum_{\m\in M_f} c_f(\m )\m$, where the symbol $\m$ denotes
a monomial,
$M_f$ is a finite set of monomials, and $c_f(\m )$ ($\m\in M_f$) are the
coefficients. For any map $\sigma$ on the field of definition
of $f$ we put $\sigma (f):=\sum_{\m\in M_f} \sigma (c_f(\m ))\m$.

We define norms for polynomials $f_i=\sum_{\m\in M_{f_i}} c_{f_i}(\m )\m$ $(i=1\kdots r )$
with complex coefficients:
\begin{eqnarray*}
&&\| f_1\kdots f_r\|:=
\max\big( |c_{f_i}(\m )|:\, 1\leq i\leq r,\, \m\in M_{f_i}\big)\, ,
\\
&&\| f_1\kdots f_r\|_1:=\sum_{i=1}^r\sum_{\m\in M_{f_i}} |c_{f_i}(\m )|
\end{eqnarray*}
and for polynomials $f_1\kdots f_r$ with coefficients in $\OQq$:
\begin{equation}\label{2.norms}
\begin{array}{l}
\| f_1\kdots f_r\|_v
:=\max \big( |c_{f_i}(\m )|_v:\, 1\leq i\leq r,\, \m\in M_{f_i}\big)
\,\,\,(v\in M_K),
\\[0.2cm]
\| f_1\kdots f_r\|_{v,1}:=
\|\sigma_v (f_1)\kdots \sigma_v(f_r)\|_1^{[K_v:\Rr ]/[K:\Qq  ]}
\,\,\,(v\in M_K^{\infty}),
\\[0.2cm]
\| f_1\kdots f_r\|_{v,1}:=
\| f_1\kdots f_r\|_v
\,\,\,(v\in M_K^0).
\end{array}
\end{equation}
Lastly, for polynomials $f_1\kdots f_r$ with coefficients in $K$ we define
heights
\begin{eqnarray*}
&&h(f_1\kdots f_r):=
\log\left( \prod_{v\in M_K}\| f_1\kdots f_r\|_v\right)\, ,
\\
&&h_1(f_1\kdots f_r):=
\log\left( \prod_{v\in M_K}\| f_1\kdots f_r\|_{v,1}\right)\, .
\end{eqnarray*}
More generally, for polynomials $f_1\kdots f_r$ with coefficients in $\OQq$
we define $h(f_1\kdots f_r)$, $h_1(f_1\kdots f_r)$
by choosing a number field $K$ containing the coefficients
of $f_1\kdots f_r$ and using the above definitions;
this is independent of the choice of $K$.

We state without proof some easy inequalities.
First, for $\x\in\OQq^{n+1}$ and
$f\in\OQq [x_0\kdots x_n]$ homogeneous
of degree $D$ we have
\begin{equation}\label{2.0}
\|f(\x )\|_v\leq \| f\|_{v,1}\| \x\|_v^D\quad\mbox{for $v\in M_K$.}
\end{equation}

Second, for $\x\in\Pp^N(\OQq )$ and
$f_0\kdots f_r\in \OQq [x_0\kdots x_N]$ homogeneous of degree $D$ we have
\begin{equation}\label{2.1*}
h({\bf y})\leq Dh({\bf x})+h_1(f_0\kdots f_r)\, ,
\end{equation}
where ${\bf y}=(f_0(\x )\kdots f_r (\x ))$.

Third, if $f\in \OQq [x_0\kdots x_n]$ is homogeneous
of degree $D$, and if $g_0\kdots g_n\in \OQq [x_0\kdots x_m]$
are homogeneous of equal degree, then for the polynomial
$f(g_0\kdots g_n)$, obtained by substituting the polynomial
$g_i(x_0\kdots x_m)$ for $x_i$
in $f$ for $i=0\kdots n$, we have
\begin{equation}
\label{2.3}
h_1\big(f(g_0\kdots g_n)\big)\leq h_1(f)+Dh_1(g_0\kdots g_n)\, .
\end{equation}

Last,
for $f_1\kdots f_r\in \OQq [x_1\kdots x_n]$ we have
\begin{equation}
\label{2.1}
h(f_1\kdots f_r)\leq h_1(f_1\kdots f_r)\leq h(f_1\kdots f_r)+\log M\, ,
\end{equation}
where $M$ is the number of non-zero coefficients in $f_1\kdots f_r$.
\end{prg}

\begin{prg}
We define another height for multihomogeneous polynomials.
Given a field $\Omega$ and tuples of non-negative integers
${\bf l}=(l_0\kdots l_m)$, we write
$\Omega [{\bf l}]$ for the set of polynomials with coefficients
in $\Omega$ in blocks of variables
${\bf z}^{(0)}=(z_0^{(0)}\kdots z_{l_0}^{(0)})$$\kdots$
${\bf z}^{(m)}=(z_0^{(m)}\kdots z_{l_m}^{(m)})$ which are
homogeneous in block ${\bf z}^{(h)}$ for $h=0\kdots m$.
For $f\in \Omega [{\bf l}]$ we denote by $\deg_h f$
the degree of $f$ in block ${\bf z}^{(h)}$.

Let
\begin{eqnarray*}
&&S(l+1):=\{ (z_0\kdots z_l)\in\Cc^{l+1}:\, |z_0|^2+\cdots +|z_l|^2=1\}\, ,
\\
&&S({\bf l}):= S(l_0+1)\times\cdots\times S(l_m+1)\, .
\end{eqnarray*}
Denote by $\mu_{l+1}$ the unique $U(l+1,\Cc )$-invariant measure
on $S(l+1)$ normalized such that $\mu_{l+1}(S(l+1))=1$,
and let $\mu_{{\bf l}}=\mu_{l_0+1}\times\cdots \times\mu_{l_m+1}$
be the product measure on $S({\bf l})$.
Then for $f\in\Cc [{\bf l}]$ we set
\begin{equation}\label{2.4}
m(f):= \int_{S({\bf l})}
\log |f({\bf z}^{(0)}\kdots {\bf z}^{(m)})|\cdot
\mu_{{\bf l}}
\,+\,\frac{1}{2}\sum_{h=0}^m \deg_h f
\left(\sum_{j=1}^{l_h}\frac{1}{2j}\right)\,.
\end{equation}
Given a number field $K$,
we define for $f\in K[{\bf l}]$,
\begin{equation}\label{2.5}
h^*(f):=\sum_{v\in M_K^{\infty}}\frac{[K_v:\Rr ]}{[K:\Qq ]}m(\sigma_v(f))
+\sum_{v\in M_K^0}\log \| f\|_{v,1}\, .
\end{equation}
Again, this does not depend on the choice of the number field $K$
containing the coefficients of $f$, so it defines a height on
$\OQq [{\bf l}]$. It is not difficult to verify that
\begin{equation}\label{2.6}
h^*(f_1\cdots f_r)=\sum_{i=1}^r h^*(f_i)
\quad\mbox{for $f_1\kdots f_r\in\OQq [{\bf l}]$.}
\end{equation}
\end{prg}
\vskip0.2cm
\begin{lemma}\label{le:2.1}
Let ${\bf l}=(l_0\kdots l_m)$ be a tuple of non-negative integers,
and $f\in\OQq [{\bf l}]$, $f\not= 0$. Then
\[
|h^*(f)-h_1(f)|\leq \sum_{h=0}^m (\deg_h f)\log (l_h+1)\, .
\]
\end{lemma}
\emph{Proof.} Put $A:=\prod_{h=0}^m (l_h+1)^{\deg_h f}$.
According to the definitions of $h^*$ and $h_1$, it suffices
to prove that for $f\in\Cc [{\bf l}]$,
\begin{equation}\label{2.7}
|m(f)-\log \| f\|_1|\leq \log A .
\end{equation}
Using $|f({\bf z}^{(0)}\kdots {\bf z}^{(m)})|\leq \| f\|_1$ for
$({\bf z}^{(0)}\kdots {\bf z}^{(m)})\in S({\bf l})$ we obtain at once
\[
m(f)\leq \log \| f\|_1+\frac{1}{2}\sum_{h=0}^m \deg_h f
\left(\sum_{j=1}^{l_h}\frac{1}{2j}\right)
\leq
\log\| f\|_1+\log A\, .
\]
To prove the inequality in the other direction, write
$f=\sum_{{\bf m}\in M_f} c({\bf m}){\bf m}$, where the sum is over a finite
number of monomials
${\bf m}=\prod_{h=0}^m\prod_{j=0}^{l_h} (z_j^{(h)})^{a_{hj}}$ with
$\sum_{j=0}^{l_h} a_{hj}=\deg_h f$ for $h=0\kdots m$. For each such monomial
we put
\[
\alpha ({\bf m}):=
\prod_{h=0}^m \frac{(\deg_h f)!}{a_{h0}!\cdots a_{h,l_h}!}\, .
\]
Then by an argument on \cite[pp. 111,112]{R2},
\[
\left( \sum_{{\bf m}\in M_f}
\alpha ({\bf m})^{-1}|c({\bf m})|^2\right)^{1/2}
\leq A^{1/2}\exp (m(f))\, .
\]
On combining this with the Cauchy-Schwarz inequality
and $\sum_{{\bf m}} \alpha ({\bf m})\leq A$, we obtain
\begin{eqnarray*}
\| f\|_1 &=&\sum_{{\bf m}\in M_f} |c({\bf m})|
\leq
\left(\sum_{{\bf m}\in M_f} \alpha ({\bf m})\right)^{1/2}
\cdot\left( \sum_{{\bf m}\in M_f}
\alpha ({\bf m})^{-1}|c({\bf m})|^2\right)^{1/2}
\\
&\leq & A\exp (m(f))\, .
\end{eqnarray*}
This proves $\log \| f\|_1\leq m(f)+\log A$, hence
\formref{2.7}.\qed
\vskip0.3cm
\begin{lemma}\label{le:2.2}
Let $f_1\kdots f_r\in\OQq [{\bf l}]$ and $f=\prod_{i=1}^r f_i$. Then
\[
h_1(f)\leq \sum_{i=1}^r h_1(f_i)\leq
h_1(f)+2\sum_{h=0}^m (\deg_h f)\log (l_h+1)\, .
\]
\end{lemma}
\emph{Proof.} The first inequality is straightforward
while the second follows from
\lemref{le:2.1} and \formref{2.6}.\qed
\vskip0.3cm
\begin{prg}
In this subsection, $X$ is a
projective subvariety of $\Pp^N$ of dimension $n\geq 1$ and degree
$d$ defined over $\OQq$.

Let $\Delta$ be a positive integer.
Denote by $M_{\Delta}$ the collection of all monomials
of degree $\Delta$ in the variables $x_0\kdots x_N$. Let
${\bf u}^{(h)}=(u_{{\bf m}}^{(h)}:\, {\bf m}\in M_{\Delta})$ $(h=0\kdots n)$
be blocks of variables.
There is an up to a constant factor unique, irreducible
polynomial $F_{X,\Delta}\in \OQq [{\bf u}^{(0)}\kdots {\bf u}^{(n)}]$,
called the $\Delta$-Chow form of $X$,
having the following property (see \cite{R1}):
\\[0.2cm]
$F_{X,\Delta}({\bf u}^{(0)}\kdots {\bf u}^{(n)})=0$ if and only if
there is a $\OQq$-rational point in the intersection of $X$ and
the hypersurfaces
$\sum_{{\bf m}\in M_{\Delta}} u_{{\bf m}}^{(h)}{\bf m}=0$ ($h=0\kdots n$).
\\[0.2cm]
Notice that $F_{X,1}$ is none other than the Chow form $F_X$ of
$X$. The form $F_{X,\Delta}$ corresponds to the Chow form
$F_{\varphi_\Delta(X)}$ of the image of $X$ under the Veronese
embedding of degree $\Delta.$ It is known that $F_{X,\Delta}$ is
homogeneous of degree $\Delta^n d$ in ${\bf u}^{(h)}$ for
$h=0\kdots n$.

For a monomial
${\bf m}=x_0^{a_0}\cdots x_N^{a_N}$ of degree $\Delta$, put
$\beta ({\bf m})=\Delta !/a_0!\cdots a_N!$. Then the modified Chow form
$G_{X,\Delta}({\bf u}^{(0)}\kdots {\bf u}^{(n)})$ is obtained by substituting
$\beta ({\bf m})^{1/2}u_{{\bf m}}^{(h)}$ for the variable $u_{{\bf m}}^{(h)}$
in the polynomial $F_{X,\Delta}({\bf u}^{(0)}\kdots {\bf u}^{(n)})$.
Notice that $G_{X,1}=F_{X,1}=F_X$. Further,
using the estimates $|\beta ({\bf m})|\leq \Delta !$,
$|\beta ({\bf m})|_p\geq |\Delta ! |_p$ for each prime number $p$,
one easily obtains
\begin{eqnarray}\label{2.8}
|h_1(F_{X,\Delta})-h_1(G_{X,\Delta})|
&\leq& \frac{1}{2}(n+1)d\Delta^n\log (\Delta !)
\\
\nonumber
&\leq& \frac{1}{2}(n+1)d\Delta^{n+1}\log \Delta \, .
\end{eqnarray}
The following is a special case of a fundamental result of R\'{e}mond
\cite[Thm. 2, pp. 99,100]{R2}:
\end{prg}
\vskip0.3cm
\begin{lemma}\label{le:2.3}
$h^*(G_{X,\Delta})=\Delta^{n+1}h^*(G_{X,1})=\Delta^{n+1}h^*(F_X)$.
\end{lemma}
\vskip0.4cm
From this we deduce:
\vskip 0.3cm
\begin{lemma}\label{le:2.4}
$h_1(F_{X,\Delta})\leq \Delta^{n+1}h(F_X)+5(n+1)d\Delta^{n+1}\log (N+\Delta )$.
\end{lemma}
\vskip 0.1cm
\emph{Proof.} Recall that $F_{X,\Delta}$ and $G_{X,\Delta}$ are
homogeneous of degree $\Delta^n d$ in each block of variables ${\bf u}^{(h)}$
$(h=0\kdots n)$ and that each of these blocks has
$\binom{N+\Delta}{\Delta} \leq (N+\Delta)^{\Delta}$ variables (that is,
the number of
coefficients of a homogeneous polynomial of degree $\Delta$ in $N+1$
variables). So by \formref{2.8} and \lemref{le:2.1},
\begin{eqnarray*}
h_1(F_{X,\Delta}) &&\leq h_1(G_{X,\Delta})
+\frac{1}{2}(n+1)d\Delta^{n+1}\log \Delta
\\
&&\leq h^*(G_{X,\Delta})+\frac{1}{2}(n+1)d\Delta^{n+1}\log \Delta
+(n+1)d\Delta^n\log\textstyle{\binom{N+\Delta}{\Delta}}
\\
&&\leq h^*(G_{X,\Delta})+\frac{3}{2}(n+1)d\Delta^{n+1}\log (N+\Delta )\, .
\end{eqnarray*}
Then using \lemref{le:2.2}, again \lemref{le:2.1} and inequality \formref{2.1}
we obtain
\begin{eqnarray*}
h_1(F_{X,\Delta}) &&\leq \Delta^{n+1}h^*(F_X)+
\frac{3}{2}(n+1)d\Delta^{n+1}\log (N+\Delta )
\\
&&\leq \Delta^{n+1}h_1(F_X)+\frac{5}{2}(n+1)d\Delta^{n+1}\log (N+\Delta )
\\
&&\leq \Delta^{n+1}h(F_X)+\frac{5}{2}(n+1)d\Delta^{n+1}\log (N+\Delta )
+\Delta^{n+1}\log M\, ,
\end{eqnarray*}
where $M$ is the number of non-zero coefficients of $F_X$.
Since $F_X$ is a polynomial
in $n+1$ blocks of $N+1$ variables,
and homogeneous of degree $d$ in each block,
we have, using \formref{5.binom}
\begin{eqnarray*}
M &&\leq \binom{N+d}{d}^{n+1}\leq \big(e(N+1)\big)^{(n+1)d}
\\
&&\leq \exp\left( \frac{5}{2}(n+1)d\log (N+\Delta )\right)\, .
\end{eqnarray*}
By inserting this into the last inequality, our lemma follows.\qed
\vskip0.4cm
We arrive at the following:
\vskip0.3cm
\begin{prop}\label{pr:2.5}
Let $g_0\kdots g_R$ be homogeneous polynomials of degree $\Delta$
in $\OQq [x_0\kdots x_N]$ such that
\[
X(\OQq )\cap \big\{g_0=0\kdots g_R=0\big\}\,=\,\emptyset\, .
\]
Let $Y=\varphi (X)$,
where $\varphi$ is the morphism on $X$ given by
$\x\mapsto (g_0(\x )\kdots g_R(\x ))$. Then
\begin{eqnarray*}
h(Y) &\leq&\Delta^{n+1}h(X)+(n+1)d\Delta^nh_1(g_0\kdots g_R)+
\\
&&\quad +5(n+1)d\Delta^{n+1}\log (N+\Delta )+3(n+1)d\Delta^n\log (R+1)\, .
\end{eqnarray*}
\end{prop}
\vskip0.15cm
\emph{Proof.}
For $j=0\kdots R$ write $y_j$ for $g_j(\x )$ and denote by ${\bf g}_j$ the
vector of coefficients of $g_j$, i.e.,
$g_j=\sum_{\m\in M_{\Delta}} c_{g_j}(\m )\m$ and
${\bf g}_j=(c_{g_j}(\m):\, \m\in M_{\Delta})$.
Introduce blocks of variables
${\bf v}^{(h)}=(v_0^{(h)}\kdots v_R^{(h)})$ $(h=0\kdots n)$ and define the
polynomial
\[
G({\bf v}^{(0)}\kdots {\bf v}^{(n)}):=
F_{X,\Delta}
\Big(\sum_{j=0}^R v_j^{(0)}{\bf g}_j\kdots
\sum_{j=0}^R v_j^{(n)}{\bf g}_j\Big)\, .
\]
Then $G({\bf v}^{(0)}\kdots {\bf v}^{(n)})=0$ if and only if
$X$ and the hypersurfaces $\sum_{j=0}^R v_j^{(h)}g_j\\=0$ ($h=0\kdots n)$
have a $\OQq$-rational point in common, if and only if
$Y$ and the hyperplanes $\sum_{j=0}^R v_j^{(h)}y_j=0$
($h=0\kdots n$) have a $\OQq$-rational
point in common, if and only if
$F_Y({\bf v}^{(0)}\kdots {\bf v}^{(n)})=0$,
where $F_Y$ is the Chow form of $Y$.
Therefore, $G$ is up to a constant factor equal to a power of $F_Y$.

Put $A:=(n+1)d\Delta^{n+1}\log (N+\Delta )$, $B:= (n+1)d\Delta^n\log (R+1)$.
Notice that $G$ has degree $d\Delta^n$ in each block ${\bf v}^{(h)}$.
Further, by \formref{2.3} we have
$h_1(G)\leq h_1(F_{X,\Delta})+(n+1)d\Delta^nh_1(g_0\kdots g_R)+B$.
Together with \lemref{le:2.2}, \lemref{le:2.1},
this implies
\begin{eqnarray*}
h(Y) &=& h(F_Y)\leq h_1(F_Y)\leq h_1(G)+2B
\\
&\leq&
h_1(F_{X,\Delta})+(n+1)d\Delta^nh_1(g_0\kdots g_R)+3B
\\
&\leq& \Delta^{n+1}h(X)+(n+1)d\Delta^nh_1(g_0\kdots g_R)+5A+3B\, ,
\end{eqnarray*}
proving our Proposition.\qed

\vskip0.5cm
\section{Proof of \thmref{th:1.2}.}\label{3}

\begin{prg}
We start with some auxiliary results.
We denote the coordinates of $\Pp^R$
by ${\bf y}=(y_0\kdots y_R)$.
\end{prg}
\vskip0.3cm
\begin{lemma}\label{le:3.1}
Let $Y$ be a projective subvariety of $\Pp^R$ of dimension $n\geq 1$
and degree $D$, defined over $\OQq$. Let ${\bf c}=(c_0\kdots c_R)$
be a tuple of reals. Let $\{ i_0\kdots i_n\}$ be a subset
of $\{ 0\kdots R\}$ such that
\begin{equation}\label{3.1}
Y(\OQq )\cap \big\{ y_{i_0}=0\kdots y_{i_n}=0\big\}\,=\,\emptyset\, .
\end{equation}
Then
\begin{equation}\label{3.2}
\frac{1}{(n+1)D}\cdot e_Y({\bf c})\geq
\frac{1}{n+1}\cdot (c_{i_0}+\cdots +c_{i_n}).
\end{equation}
\end{lemma}
\vskip0.2cm
\emph{Proof.}
For a subset $I=\{k_0\kdots k_n\}$ of $\{0\kdots R\}$ with
$k_0<k_1<\cdots<k_n,$ define the {\it bracket}
\[
[I]=[I]({\bf u}^{(0)}\kdots {\bf u}^{(n)})
:=\operatorname{det}\left(u_{k_j}^{(i)}\right)_{i,j=0\kdots n},
\]
where again ${\bf u}^{(h)}$ denotes the block of variables
$(u_0^{(h)}\kdots u_R^{(h)})$.
Let $I_1\kdots I_S$ with $S=\binom{R+1}{n+1}$ be all subsets of
$\{0\kdots R\}$ of cardinality $n+1.$ Then the Chow form $F_Y$
of $Y$ can be written
as a homogeneous polynomial of degree $D$ in $[I_1]\kdots [I_S]:$
\begin{eqnarray}\label{3.3}
F_Y=\sum_{{\bf a}\in A}C({\bf a})[I_1]^{a_1}
\cdots [I_S]^{a_S},
\end{eqnarray}
where $A$ is the set of tuples of non-negative
integers ${\bf a}=(a_1\kdots a_S)$ with $a_1+\cdots+a_S=D$
and where $C({\bf a})\in\OQq$ for ${\bf a}\in A$
\cite[p. 41, Theorem IV]{HP}.
For each bracket $[I]$ we have
\[
[I](t^{c_0}u^{(0)}_0\kdots t^{c_R}u^{(0)}_R;\ldots;
t^{c_0}u^{(n)}_0
\kdots t^{c_R}u^{(n)}_R)=
t^{\sum_{i\in I}c_i}[I],
\]
therefore,
\begin{eqnarray}\label{3.3a}
&&F_Y(t^{c_0}u^{(0)}_0\kdots t^{c_R}u^{(0)}_R;\ldots;
t^{c_0}u^{(n)}_0
\kdots t^{c_R}u^{(n)}_R)\\
\nonumber
&&
\qquad\qquad\qquad
=\sum_{{\bf a}\in A}C({\bf a})t^{\sum_{j=1}^Sa_j(\sum_{i\in I_j} c_i)}[I_1]^{a_1}
\cdots [I_S]^{a_S}\, .
\end{eqnarray}
Put ${\bf e}_0:=(1,0\kdots 0)$, ${\bf e}_1:=(0,1\kdots 0)$,
$\ldots$,
${\bf e}_R:=(0,0\kdots 1)$. Write
$\{ i_0\kdots i_n\}=: I_1$.
By \formref{3.1} we have
$F_Y({\bf e}_{i_0}\kdots {\bf e}_{i_n})\neq 0.$ Further,
\[
[I_1]({\bf e}_{i_0}\kdots {\bf e}_{i_n})=1,\quad
[I]({\bf e}_{i_0}\kdots {\bf e}_{i_n})=0\text{ for }I\neq I_1.
\]
Hence
in expression \formref{3.3} there is a
term $C\cdot [I_1]^D$ with $C\in\OQq^*$,
and if we substitute ${\bf u}^{(j)}={\bf e}_{i_j}$ $(j=0\kdots n)$
in \formref{3.3a} we obtain $C\cdot t^{D(c_{i_0}+\cdots +c_{i_n})}$.
That is,
one of the numbers $e_i$ in \formref{4.4} is equal to
$D\left(c_{i_0}+\cdots +c_{i_n}\right)$.
This implies \formref{3.2} at once. \qed
\vskip0.4cm

In addition, we need the following combinatorial lemma, which is
a consequence of \cite[Lemma 4]{E1}.
\vskip0.3cm
\begin{lemma}\label{le:3.3}
Let $\theta$ be a real with $0<\theta\leq\frac{1}{2}$ and let $q$ be a positive
integer. Then there exists a set $\WW$ of cardinality at most
$(e/\theta )^{q-1}$, consisting of tuples $(c_1\kdots c_q)$ of non-negative
reals with $c_1+\cdots +c_q=1$, with the following property:\\
for every set of reals $A_1\kdots A_q$ and $\Lambda$ with
$A_j\leq 0$ for $j=1\kdots q$ and $\sum_{j=1}^q A_j\leq -\Lambda$, there exists
a tuple $(c_1\kdots c_q)\in\WW$ such that
\[
A_j\leq -c_j(1-\theta )\Lambda\,\,\,\mbox{for $j=1\kdots q$.}
\]
\end{lemma}
\vskip0.4cm
\begin{prg}
In what follows,
$K$ is a number field, $S$ a finite set of places of $K$,
and $X$, $N$, $n$, $d$, $s$, $C$,
$f_i^{(v)}$ ($v\in S$, $i=0\kdots n$), $C$, $\Delta$, $A_1$, $A_2$, $A_3$, $H$
are as in \thmref{th:1.2}. We denote the coordinates on $\Pp^N$ by
${\bf x}=(x_0\kdots x_N)$.

Let $f_0\kdots f_R$ be the distinct polynomials among
$\sigma (f_j^{(v)})$ ($v\in S$, $j=0\kdots n$, $\sigma\in G_K$).
Then by \formref{1.5a},
\begin{equation}\label{3.4}
R\leq C(n+1)s-1\, .
\end{equation}
Let $K'$ be the extension of $K$ generated by the coefficients
of $f_0\kdots f_R$.
Put $g_i:= f_i^{\Delta /\deg f_i}$ for $i=0\kdots R$.
Thus, $g_0\kdots g_R$ are homogenenous polynomials
in $K' [x_0\kdots x_N]$ of degree $\Delta$. Define
\[
\varphi :\, {\bf x}\mapsto (g_0({\bf x})\kdots g_R({\bf x})),\quad
Y:=\varphi (X)\, .
\]
By assumption \formref{1.3},
$\varphi$ is a
finite morphism on $X$, and $Y$ is a projective subvariety
of $\Pp^R$ defined over $K'$.
We have
\begin{equation}\label{3.5}
\dim Y=n,\quad \deg Y=: D\leq d\Delta^n\, .
\end{equation}

We denote places on $K'$ by $v'$
and define normalized absolute values $|\cdot |_{v'}$ on $K'$ similarly
to \prgref{notation}. Further, for every $v'\in M_{K'}$ we choose an
extension of $|\cdot |_{v'}$ to $\OQq$.
Since $K'/K$
is a normal extension, for every $v'\in M_{K'}$
there is $\tau_{v'}\in G_K$ such that
\begin{equation}\label{3.5a}
|x|_{v'}=|\tau_{v'} (x)|_v^{1/g(v)}\quad\mbox{for $x\in \OQq$}
\end{equation}
where $v\in M_K$ is the place below $v'$ and $g(v)$ is the number of places
of $K'$ lying above $v$.
For each $v'\in M_{K'}^{\infty}$ there is an isomorphic
embedding $\sigma_{v'}:\, K'\hookrightarrow \Cc$ such that
$|x|_{v'}=|\sigma_{v'}(x)|^{[K_{v'}:\Rr ]/[K':\Qq ]}$ for $x\in\OQq$.
We define norms $\|\cdot \|_{v'}$, $\|\cdot \|_{v' ,1}$ for polynomials
similarly as in \formref{2.norms}, with $K',v',\sigma_{v'}$ in place of
$K,v,\sigma_v$.
\end{prg}

\begin{prg}
For later purposes we estimate from above $h_1(1,g_0\kdots g_R)$ and $h(Y)$.
By a straightforward computation we have
for $v'\in M_{K'}^{\infty}$,
\begin{eqnarray*}
&&\| 1, \sigma_{v'} (g_0)\kdots \sigma_{v'}(g_R)\|_1
\\
&&\qquad =1+\sum_{i=0}^R \| \sigma_{v'} (g_i)\|_1
\leq 1+\sum_{i=0}^R \| \sigma_{v'}(f_i)\|_1^{\Delta /\deg f_i}
\\
&&\qquad\leq 1+\sum_{i=0}^R
\left(\binom{\deg f_i +N}{\deg f_i}\| \sigma_{v'}(f_i)\|\right)^
{\Delta /\deg f_i}
\\
&&\qquad\leq (R+2)(N+\Delta )^{\Delta}
\| 1,\sigma_{v'}(f_0)\kdots \sigma_{v'} (f_R)\|^{\Delta}
\\
&&\qquad\leq (R+2)(N+\Delta )^{\Delta}
\prod_{i=0}^R\| 1,\sigma_{v'}(f_i)\|^{\Delta}\, .
\end{eqnarray*}
So for $v'\in M_{K'}^{\infty}$ we have
\[
\| 1,g_0\kdots g_R\|_{v',1}\leq \left( (R+2)(N+\Delta )^{\Delta}\right)
^{\frac{[K'_{v'}:\Rr ]}{[K':\Qq]}}\cdot
\prod_{i=0}^R \| 1,f_i\|_{v'}^{\Delta}.
\]
In an easier manner one obtains for $v'\in M_{K'}^0$,
\[
\| 1,g_0\kdots g_R\|_{v',1}\leq \prod_{i=0}^R\| 1,f_i\|_{v'}^{\Delta}.
\]
So by taking the product over $v'\in M_{K'}$, substituting \formref{3.4},
and using that polynomials with conjugate sets of coefficients have the
same height,
\begin{eqnarray*}
&&
h_1(1,g_0\kdots g_R)\leq \Delta \big(\sum_{i=0}^R h(1,f_i)\big)+
\Delta\log \big((R+2)(N+\Delta )^{\Delta}\big)
\\
&&
\quad\quad\leq \Delta C\big(\sum_{v\in S}\sum_{j=0}^n h(1,f_j^{(v)})\big)+
\Delta\log (N+\Delta )+
\log (3Cns)\, ,
\end{eqnarray*}
and by inserting this estimate into \propref{pr:2.5} we infer
\begin{eqnarray*}
&&h(Y)\leq\Delta^{n+1}h(X)+(n+1)d\Delta^{n+1}C\sum_{v\in S}\sum_{j=0}^n h(1,f_j^{(v)})+
\\
&&\qquad\qquad +\, 6(n+1)d\Delta^{n+1}\log (N+\Delta )+4(n+1)d\Delta^n\log (3Cns)\, .
\end{eqnarray*}
A straightforward computation gives the more tractable estimates
\begin{eqnarray}\label{3.6}
&&h_1(g_0\kdots g_R)\leq 6\Delta^2Cns\cdot H\, ,
\\
\label{3.7}
&& h(Y)\leq 25n^2d\Delta^{n+2}Cs\cdot H\, ,
\end{eqnarray}
where $H$ is defined by \formref{1.5b}.
\end{prg}

\begin{prg}
We reduce \formref{1.6} to a finite number of systems of inequalities,
and then show that each such system leads to an inequality involving
a twisted height.

Let $\x\in X(\OQq )$ be a solution of \formref{1.6}.
For $v\in S$, let $I_v$ be the subset of $\{ 0\kdots R\}$
such that $\{ f_j^{(v)}:\, j=0\kdots n\}=\{ f_i:\, i\in I_v\}$.
Put
$G_v:=\| 1,g_0\kdots g_R\|_{v,1}$ for $v\in S$. Then
\[
\sum_{v\in S}\sum_{i\in I_v} \log\left( \max_{\sigma\in G_K}
\frac{|g_i(\x )|_v}{G_v\|\sigma (\x )\|_v^{\Delta}}\right) \leq
-(n+1+\delta )\Delta h(\x )\, .
\]
By \formref{2.0}, the terms in the sum are $\leq 0$.
We apply \lemref{le:3.3}
with $q=(n+1)s$ and
$\theta =\frac{\delta}{(2n+2+2\delta )}=1-\frac{n+1+\delta /2}{n+1+\delta}$.
We infer that there is a set
$\WW$ with
\begin{equation}\label{3.10}
\#\WW \leq \left(\frac{e(2n+2+2\delta )}{\delta}\right)^{(n+1)s-1}
\leq (17n\delta^{-1})^{(n+1)s-1}
\end{equation}
consisting of tuples of non-negative reals $(c_{iv}:\, v\in S,\, i\in I_v)$
with
\begin{equation}\label{3.11a}
\sum_{v\in S}\sum_{i\in I_v} c_{iv}=1\, ,
\end{equation}
such that for every solution ${\bf x}\in X(\OQq )$ of \formref{1.6}
there is a tuple $(c_{iv}:\, v\in S,\, i\in I_v)\in\WW$ with
\begin{eqnarray}\label{3.12a}
\log\left(\max_{\sigma\in G_K} \frac{|g_i(\sigma (\x
))|_v}{G_v\cdot \|\sigma (\x )\|_v^{\Delta}}\right) &\leq&
-c_{iv}\left(n+1+\frac{\delta}{2}\right)\Delta h(\x )
\\
\nonumber
&&\quad\quad\quad\quad\quad
\mbox{($v\in S$, $i\in I_v$).}
\end{eqnarray}

Denote by $S'$ the set of places of $K'$ lying above the places in
$S$. Notice that each element of $G_K$ acts as a permutation on
$g_0\kdots g_R$. Let $v'\in S'$.
Write $v$ for the place of $K$
lying below $v'$ and let $\tau_{v'}\in G_K$
be given by \formref{3.5a}.
Then we define $I_{v'}\subset\{ 0\kdots R\}$,
$c_{i,v'}$ $(i\in I_{v'})$ by
\begin{eqnarray*}
&&
\{ g_i :\, i\in I_{v'}\}= \{ \tau_{v'}^{-1}(g_j):\, j\in I_v\}
\quad\mbox{for $v'\in S'$,}
\\
&&
c_{i,v'}:= c_{jv}/g(v)\quad\mbox{for $v'\in S'$, $i\in I_{v'}$,}
\end{eqnarray*}
where $j\in I_v$ is the index such that
$g_i=\tau_{v'}^{-1}(g_j)$.
Further, we put
\[
G_{v'}:=\| 1,g_0\kdots g_R\|_{v',1}\quad\mbox{for $v'\in M_{K'}$.}
\]
Then in view of \formref{3.5a}, we can rewrite system \formref{3.12a} as
\begin{eqnarray}\label{3.12}
\log\left(\max_{\sigma\in G_K} \frac{|g_i(\sigma (\x
))|_{v'}}{G_{v'}\cdot \|\sigma (\x )\|_{v'}^{\Delta}}\right)
&\leq& -c_{i,v'}\left(n+1+\frac{\delta}{2}\right)\Delta h(\x )
\\
\nonumber
&&\quad\quad\quad\quad\quad
\mbox{($v'\in S'$, $i\in I_{v'}$).}
\end{eqnarray}
Invoking \formref{3.10}, \formref{3.11a} we obtain the following:
\end{prg}
\vskip0.3cm
\begin{lemma}\label{le:3.2}
There is a set $\WW '$ of cardinality at most $(17n\delta^{-1})^{(n+1)s-1}$,
consisting of tuples of non-negative reals
$(c_{i,v'}:\, v'\in S',\,i\in I_{v'})$ with
\begin{equation}\label{3.11}
\sum_{v\in S'}\sum_{i\in I_{v'}} c_{i,v'}=1\, ,
\end{equation}
with the property that for every $\x\in X(\OQq )$ with \formref{1.6} there is a
tuple in $\WW '$ such that $\x$ satisfies \formref{3.12}.
\end{lemma}
\vskip0.2cm

We consider the solutions of a fixed system \formref{3.12}.
Put
\begin{eqnarray}\label{3.11b}
&c_{i,v'}=0\,\,&\mbox{for $v'\in S'$, $i\in\{ 0\kdots R\}\backslash I_{v'}$}
\\
\nonumber
&&\mbox{and $v'\in M_{K'}\backslash S'$, $i=0\kdots R$}
\end{eqnarray}
and put ${\bf c}_{v'}:=(c_{0,v'}\kdots c_{R,v'})$ for $v'\in M_{K'}$,
${\bf c}:= ({\bf c}_{v'}:\, v'\in M_{K'})$.
Denote by ${\bf y}=(y_0\kdots y_R)$ the coordinates of $\Pp^R$.
We define $H_{Q,{\bf c}}({\bf y})$,
$E_Y({\bf c})$ similarly as
\formref{4.3}, \formref{4.6a}, respectively,
but with $K'$ in place of $K$.
\vskip0.3cm
\begin{lemma}\label{le:3.4}
Let ${\bf x}\in X(\OQq )$ be a solution of \formref{3.12} satisfying
\formref{1.7} and let $\sigma\in G_K$.
Put
\[
{\bf y}:=\varphi (\sigma ({\bf x})),\quad
Q:=\exp\Big( (n+1+\delta/2)\Delta h({\bf x})\Big).
\]
Then
\begin{equation}\label{3.14}
H_{Q,{\bf c}}({\bf y})\leq Q^{E_Y({\bf c})-\frac{\delta}{2(n+2)^2}}\, .
\end{equation}
\end{lemma}
\vskip0.2cm
\emph{Proof.} We first estimate from below $E_Y({\bf c})$.
Let $v'\in S'$ and write $I_{v'}=\{ i_0\kdots i_n\}$.
From assumption \formref{1.3},
and from the fact that $X$ is defined over $K$ and that
$g_{i_0}\kdots g_{i_n}$ are conjugate
over $K$ to powers of $f_0^{(v)}\kdots f_n^{(v)}$ where $v\in S$ is the
place below $v'$, it follows that
$X(\OQq )\cap \big\{ g_{i_0}=0\kdots g_{i_n}=0\big\}\,=\emptyset$.
Since $Y=\varphi (X)$, for ${\bf y}\in Y(\OQq )$
there is $\x\in X(\OQq )$ with $y_i=g_i(\x )$
for $i=0\kdots R$. Hence
\[
Y(\OQq )\cap \big\{ y_{i_0}=0\kdots y_{i_n}=0\big\}\,=\,\emptyset\, .
\]
Now \lemref{le:3.1} implies
\[
\frac{1}{(n+1)D}\cdot e_Y({\bf c}_{v'})\geq
\frac{1}{n+1}(c_{i_0,v'}+\cdots +c_{i_n,v'})
=\frac{1}{n+1}\cdot\sum_{i\in I_{v'}} c_{i,v'}\, .
\]
This holds for $v'\in S'$. For $v'\not\in S'$ we have
$e_Y({\bf c}_{v'})=0$ by \formref{3.11b}.
By summing over $v'\in S'$ and using \formref{3.11}, we arrive at
\begin{equation}\label{3.13}
E_Y({\bf c})\geq \frac{1}{n+1}\, .
\end{equation}

Now let $\x\in X(\OQq )$ be a solution of \formref{3.12} with \formref{1.7}
and let $\sigma\in G_K$.
Then $\sigma (\x )$ is also a solution of \formref{3.12}. In fact,
by \formref{3.11b},
$\sigma (\x )$ satisfies \formref{3.12}
for $v\in M_K$, $i=0\kdots R$. Write ${\bf y}=\varphi (\sigma (\x ))$
so that $y_i=g_i(\sigma (\x ))$ for $i=0\kdots R$.
Let $L$ be a finite normal extension of $K'$ such that
$\sigma (\x )\in X(L)$. Pick $w\in M_L$ and let $v'$ be the place of $K'$
below $w$.
Then there is $\tau_w\in {\rm Gal}(\OQq /K')$ such that
$|x|_w=|\tau_w(x)|_{v'}^{d(w|v')}$ for $x\in L$,
where $d(w|v')=[L_w:K'_{v'}]/[L:K' ]$.
Hence
for $i=0\kdots R$, with the usual notation $c_{iw}=d(w|v')c_{i,v'}$,
\begin{eqnarray*}
|y_i|_wQ^{c_{iw}}&=&|g_i(\sigma (\x ))|_wQ^{c_{iw}}
=\left(|g_i(\tau_w\sigma ({\bf x}))|_{v'}Q^{c_{i,{v'}}}\right)^{d(w|v')}
\\
&\leq&\left( G_{v'}\|\tau_w\sigma ({\bf x})\|_{v'}^{\Delta}\right)^{d(w|v')}
=G_{v'}^{d(w|v)}\|\sigma ({\bf x})\|_w^{\Delta}\, .
\end{eqnarray*}
By taking the product over $w\in M_L$ and using $h(\sigma (\x ))= h(\x )$
we obtain
\[
H_{Q,{\bf c}}({\bf y})\leq \exp (h_1(1,g_0\kdots g_R))\cdot
Q^{\frac{1}{n+1+\delta /2}}\, .
\]
Now \formref{3.14} follows by observing that by \formref{3.13},
assumption \formref{1.7},
and \formref{3.6},
\begin{eqnarray*}
&&\left (E_Y({\bf c})-\frac{\delta}{2(n+2)^2}-\frac{1}{n+1+\delta /2}\right)
\log Q
\\
&&\qquad\geq
\left(\frac{1}{n+1}-\frac{\delta}{2(n+2)^2}-\frac{1}{n+1+\delta /2}\right)\log Q
\\
&&\qquad
=\frac{\delta (4n+6-\delta (n+1))}{4(n+1)(n+2)^2}\cdot\Delta h(\x )
\geq\frac{\delta\Delta}{2(n+2)^2}A_3H
\\
&&\qquad\geq 6\Delta^2CnsH\geq h_1(1,g_0\kdots g_R)\, .
\end{eqnarray*}
\qed
\vskip0.2cm
\begin{prg}
We finish the proof of \thmref{th:1.2}.
We apply \thmref{th:4.1} with $K'$, $\frac{\delta}{2(n+2)^2}$ in place of
$K$, $\delta$ and, in view of \formref{3.4} and \formref{3.5}, with
$D\leq d\Delta^n$ and $R=C(n+1)s-1$.
Notice that by \formref{3.11},\formref{3.11b}, the conditions
\formref{4.6}, \formref{4.6b}, \formref{4.6c} (with $K'$ in place of $K$)
are satisfied.
Denote by $B_1'$, $B_2'$, $B_3'$
the quantities obtained by substituting
$\frac{\delta}{2(n+2)^2}$ for $\delta$,
$C(n+1)s-1$ for $R$, and $d\Delta^n$ for $D$ in
the quantities $B_1$, $B_2$, $B_3$,
respectively, defined by \formref{4.7}.
Recall that if ${\bf x}$ satisfies \formref{1.7}
then \lemref{le:3.4} is applicable. Moreover,
\begin{eqnarray*}
\log Q &=& \left(n+1+\frac{\delta}{2}\right)\Delta h({\bf x}) \geq
A_3H
\\
&=& \exp\big( 2^{6n+20}n^{2n+3}\delta^{-n-1}d^{n+2}\Delta^{n(n+2)}\log (2Cs)\big)\cdot H
\\[0.2cm]
&\geq& \exp\left( 2^{5n+4}(2(n+2)^2\delta^{-1})^{n+1}(d\Delta^n)^{n+2}
\log (4C(n+1)s)\right)\cdot
\\
&&\quad
\cdot \big( 26n^2d\Delta^{n+2}Cs\big)\cdot H
\\[0.2cm]
&=& B_3'\cdot \big( 26n^2d\Delta^{n+2}Cs\big)\cdot H\geq B_3'(h(Y)+1)\, ,
\end{eqnarray*}
where the last inequality follows from \formref{3.7}.
Hence \thmref{th:4.1} is applicable.

Now \thmref{th:4.1} and \lemref{le:3.4} imply that
there are homogeneous polynomials $F_1\kdots F_t\in K'[y_0\kdots y_R]$
not vanishing identically on $Y$,
with $t\leq B_1'$ and $\deg F_i\leq B_2'$ for $i=1\kdots t$,
with the property that
for every solution ${\bf x}\in X(\OQq )$ of \formref{3.12} with
\formref{1.7}, there is $F_i\in\{ F_1\kdots F_t\}$
such that
$F_i(\varphi (\sigma (\x )))=0$ for every $\sigma\in G_K$.
(In fact, taking $Q=\exp\big( (n+1+\delta /2 )\Delta h(\x )\big)$
it follows from \thmref{th:4.1} that there is $F_i$
with $F_i({\bf y})=0$ for every ${\bf y}\in Y(\OQq )$
with $H_{Q,{\bf c}}({\bf y})\leq Q^{E_Y({\bf c})-\delta /2(n+2)^2}$,
and then by \lemref{le:3.4} this holds in particular for
all points ${\bf y}=\varphi (\sigma (\x ))$, $\sigma\in G_K$.)

This means that $\tilde{F}_i(\sigma (\x ))=0$ for $\sigma\in G_K$,
where $\tilde{F}_i$ is the polynomial obtained by
substituting $g_j$ for $y_j$ in $F_i$
for $j=0\kdots R$. Notice that $\tilde{F}_i\in K' [x_0\kdots x_N]$,
$\deg \tilde{F}_i\leq B_2'\Delta$,
and that $\tilde{F}_i$ does not vanish identically
on $X$. Write $\tilde{F}_i=\sum_{k=1}^M \omega_k\tilde{F}_{ik}$
where $\omega_1\kdots\omega_M$ is a $K$-basis of $K'$,
and the $\tilde{F}_{ik}$ are polynomials with coefficients in $K$.
We can choose $G_i\in\{ \tilde{F}_{ik}:\, k=1\kdots M\}$ which does not vanish
identically on $X$.
Now $\sigma(\tilde{F}_i)(\x )=0$ for $\sigma\in G_K$. Since the polynomials
$\tilde{F}_{ik}$ are linear combinations of the polynomials
$\sigma (\tilde{F}_i)$ ($\sigma\in G_K$) it follows that $\tilde{F}_{ik}(\x)=0$
for $k=1\kdots M$, so in particular $G_i(\x )=0$.

It follows that there are homogeneous polynomials
$G_1\kdots G_t\in \\K[x_0\kdots x_N]$ with
$t\leq B_1'$ and $\deg G_i \leq B_2'\Delta$ for $i=1\kdots t$,
not vanishing identically on $X$, such that
the set of $\x\in X(\OQq )$ with \formref{3.12} and with \formref{1.7}
is contained in $\bigcup_{i=1}^t\Big( X\cap \{G_i=0\}\Big)$.

According to \lemref{le:3.2},
there are at most $T:= \left( 17n\delta^{-1}\right)^{(n+1)s-1}$
different systems \formref{3.12}, such that
every solution ${\bf x}\in X(\OQq )$
of \formref{1.6} satisfies one of these systems.
Consequently, there
are homogeneous polynomials $G_1\kdots G_u\in K[x_0\kdots x_N]$
not vanishing identically on $X$, with $u\leq B_1'T$ and with
$\deg G_i\leq B_2'\Delta$ for $i=1\kdots u$, such that the set of
${\bf x}\in X(\OQq )$ with \formref{1.6}, \formref{1.7} is contained
in $\bigcup_{i=1}^u\Big( X\cap \{G_i=0\}\Big)$.

Now the proof of \thmref{th:1.2}
is completed by observing that in view of \formref{4.7},
\[
B_2'\Delta = (4n+3)(d\Delta^n)(2(n+2)^2\delta^{-1})\Delta
=(8n+6)(n+2)^2d\Delta^{n+1}\delta^{-1}=A_2
\]
and
\begin{eqnarray*}
B_1'T&\leq&
\exp\left(  2^{10n+4}(2(n+2)^2)^{2n}\delta^{-2n}(d\Delta^n)^{2n+2}\right)
\cdot
\\
&&\quad
\cdot\log (4(n+1)Cs)\log\log (4(n+1)Cs)\cdot
\left( 17n\delta^{-1}\right)^{(n+1)s-1}
\\[0.2cm]
&\leq&
\exp\left( 2^{12n+16}n^{4n}\delta^{-2n}d^{2n+2}\Delta^{n(2n+2)}\right)
\cdot
\\
&&\quad
\cdot (20n\delta^{-1})^{(n+1)s}\cdot
\log (4C)\log\log (4C)
\\[0.2cm]
&=&A_1\, .
\qquad\qquad\qquad\qquad\qquad\qquad\qquad\qquad\qquad\qquad\qquad\qquad\qed
\end{eqnarray*}
\end{prg}

\end{document}